\documentclass[12pt,draftcls,onecolumn]{IEEEtran}
\usepackage{texdraw,amsmath,amssymb,bm,cite,graphicx}

\newcommand{\beq}{\begin{equation}}
\newcommand{\eeq}{\end{equation}}
\newcommand{\bea}{\begin{eqnarray}}
\newcommand{\eea}{\end{eqnarray}}
\newcommand{\bef}{\begin{figure}}
\newcommand{\eef}{\end{figure}}
\newcommand{\bd}{\begin{displaymath}}
\newcommand{\ed}{\end{displaymath}}
\newcommand{\nn}{\nonumber}
\newcommand{\nnl}{\nonumber \\}
\newcommand{\fig}[1]{Fig.\ \ref{#1}}

\newcommand{\bsc}{\begin{scriptstyle}}
\newcommand{\esc}{\end{scriptstyle}}
\def\bmat{\left[ \begin{array}}
\def\emat{\end{array} \right]}
\newcommand{\defn}{\stackrel{\triangle}{=}}

\begin{document}
\bibliographystyle{ieeetr}

\title{Robust State Space Filtering under Incremental Model 
Perturbations subject to a Relative Entropy Tolerance}

\author{Bernard C. Levy\thanks{B. Levy is with the Department
of Electrical and Computer Engineering, 1 Shields Avenue, University
of California, Davis, CA 95616. Email: bclevy@ucdavis.edu, phone:
(530) 752-8025, fax: (530) 752-8428.} ~ and ~ Ramine 
Nikoukhah\thanks{R. Nikoukhah is with the Institut de Recherche
en Informatique et Automatique (INRIA), Domaine de Voluceau,
Rocquencourt, 78153 Le Chesnay Cedex France. email: 
ramine.nikoukhah@inria.fr.}}

\date{\today}

\maketitle

\begin{abstract}

This paper considers robust filtering for a nominal Gaussian
state-space model, when a relative entropy tolerance is
applied to each time increment of a dynamical model. The problem
is formulated as a dynamic minimax game where the maximizer 
adopts a myopic strategy. This game is shown to admit a 
saddle point whose structure is characterized by applying 
and extending results presented earlier in \cite{LN} for
static least-squares estimation. The resulting minimax filter
takes the form of a risk-sensitive filter with a time varying
risk sensitivity parameter, which depends on the tolerance bound
applied to the model dynamics and observations at the corresponding 
time index. The least-favorable model is constructed and used to 
evaluate the performance of alternative filters. Simulations comparing 
the proposed risk-sensitive filter to a standard Kalman filter show a 
significant performance advantage when applied to the least-favorable 
model, and only a small performance loss for the nominal model. 

\end{abstract}

\begin{keywords}

commitment, dynamic minimax game, least-favorable model,
myopic strategy, relative entropy, risk-sensitive filtering, 
robust filtering.

\end{keywords}

\section{Introduction}
\label{sec:intro}

Soon after the introduction of Wiener and Kalman filters, it was 
recognized that these filters were vulnerable to modelling errors,
in the form of either parasitic signals or perturbations of the system
dynamics. Various approaches were proposed over the last 35 years 
to construct filters with a guaranteed level of immunity to modelling 
uncertainties. Drawing from the framework developed by Huber
for robust statistics \cite{Hub}, Kassam, Poor and their
collaborators proposed an approach \cite{KL,Poo,KP} where the
optimum filter is selected by solving a minimax problem. In 
this approach, the set of possible system models is described by
a neighborhood centered about the nominal model, and two players  
affront each other. One player (say, nature) selects the 
least-favorable model in the allowable neighborhood and the other
player designs the optimum filter for the least-favorable model. 
While minimax filtering is conceptually simple, its implementation
can be very difficult, since it depends on the specification
of the allowable neighborhood and of the loss function to be
minimized. After some early success in the design of Wiener filters
for neighborhoods specified by $\epsilon$-contamination models or
power spectral bands, progress stalled gradually and researchers 
started looking in different directions to develop robust filters. 
The 1980s saw the development of an entirely different class of
robust filters based on the minimization of risk-sensitive and
$H^{\infty}$ performance criteria \cite{SDJ,Whi,NK,MG,HSK}. This
approach seeks to avoid large errors, even if these errors are
unlikely based on the nominal model. For example, risk-sensitive
filters replace the standard quadratic loss function of least-squares
filtering by an exponential quadratic function, which of course 
penalizes severely large errors. However, errors in the model 
dynamics are not introduced explicitly in $H^{\infty}$ and 
risk-sensitive filtering, and the growing awareness of the importance
of such errors prompted a number of researchers in the early 2000s
\cite{PS,Say,GC} to revive the minimax filtering viewpoint, but
in a context where modelling errors are described in terms of 
norms for state-space dynamics perturbations. The present paper,
which is a continuation of \cite{LN}, can be viewed as part of a
larger effort initiated by Hansen and Sargent \cite{HS,HS1,HS2} and
other researchers \cite{BJP,YUP} which is aimed at reinterpreting
risk-sensitive filtering from a minimax viewpoint. In this context,
modelling uncertainties are described by specifying a tolerance for
the relative entropy between the actual system and the nominal 
model. To a fixed tolerance level describing the modeller's
confidence in the nominal model corresponds a ball of possible
models for which it is then possible to apply the minimax
filtering approach proposed by Kassam and Poor.   

The minimax formulation of robust filtering based on a relative
entropy constraint has several attractive features. First, 
relative entropy is a natural measure of model mismatch
which is commonly used by statisticans for fitting statistical
models by using techniques such as the expectation maximization
iteration \cite{MK}. More fundamentally, it was shown by Chentsov
\cite{Cen} and by Amari \cite{AN} that manifolds of statistical
models can be endowed with a non-Riemannian differential  
geometric structure involving two dual connections associated to
the relative entropy and the reverse relative entropy. In addition to
this strong theoretical justification, it turns out that minimax Wiener 
and Kalman filtering problems with a relative entropy constraint admit 
solutions \cite{YUP,LN,HS1,HS} in the form of risk-sensitive filters,
thus providing a new interpretation for such filters. The main
difference between earlier works and the present paper is that,
instead of placing a single relative entropy constraint on the entire
system model, we apply a separate constraint to each time increment
of the model. This approach, which is closer to the one advocated in
\cite{Say,GC} is based on the following consideration. When a single
relative entropy constraint is placed on the complete system model, 
the maximizing player has the opportunity to allocate almost all of
its mismatch budget to a single element of the model most susceptible
to uncertainties. But this strategy may lead to overly pessimistic
conclusions, since in practice modellers allocate the same level
of effort to modelling each time component of the system. Thus it 
probably makes better sense to specify a fixed uncertainty tolerance 
for each model increment, instead of a single bound for the overall
model.

The analysis presented relies in part on applying and extending
static least-squares estimation results derived in \cite{LN} for
nominal Gaussian models. These results are reviewed in Section
\ref{sec:static}. In Section \ref{sec:game}, the robust state-space
filtering problem with an incremental relative entropy constraint
is formulated as a dynamic game where the maximizer adopts a
myopic strategy whose goal is to maximize the mean-square 
estimation error at the current time only. The existence of a saddle 
point is established in Section \ref{sec:filter}, where the 
least-favorable model specifying the saddle point is characterized
by extending a Lemma of \cite{LN} to the dynamic case. A careful 
examination of the least-favorable model structure allows the 
transformation of the dynamic estimation game into an equivalent static 
one, to which the result of \cite{LN} become applicable. The robust 
filter that we obtain is a risk-sensitive filter, but with a time-varying
risk-sensitivity parameter, representing the inverse of the Lagrange
multiplier associated to the model component constraint for the corresponding 
time period. The least-favorable state-space model for which the 
robust filter is optimal is constructed in Section \ref{sec:leastfav}.
This model extends to the finite-horizon time-varying case a 
model derived asymptotically in \cite{HS} for the case of constant 
systems. The least-favorable state-space model allows performance 
evaluation studies comparing the performance of the minimax
filter with that of other filters, such as the standard Kalman 
filter. Simulations are presented in Section \ref{sec:simu} which 
illustrate the dependence of the filter performance on the 
relative entropy tolerance applied to each model time component.
The robust filter is compared to the ordinary Kalman filter by
examining their respective performances for both the nominal and 
least-favorable systems. Finally, some conclusions are presented in 
Section \ref{sec:conc}.

\section{Robust Static Estimation}
\label{sec:static}

We start by reviewing a robust static estimation result 
derived in \cite{LN}, since its extension to the dynamic
case is the basis for the robust filtering scheme we
propose. Let 
\beq
z = \bmat {c} 
x \\ y
\emat
\label{2.1}
\eeq
be a random vector of ${\mathbb R}^{n+p}$, where
$x \in {\mathbb R}^n$ is a vector to be estimated,
and $y \in {\mathbb R}^p$ is an observed vector. The
nominal and actual probability densities of $z$
are denoted respectively as $f (z)$ and $\tilde{f} (z)$.
The deviation of $\tilde{f}$ from $f$ is measured by
the relative entropy (Kullback-Leibler divergence)
\beq
D (\tilde{f},f) = \int_{{\mathbb R}^{n+p}} 
\tilde{f}(z) \ln \Big( \frac{\tilde{f}(z)}{f(z)} \Big)
dz \: . \label{2.2}
\eeq 
The relative entropy is not a distance, since it is
not symmetric and does not satisfy the triangle inequality.
However it has the property that $D(\tilde{f},f) \geq 0$
with equality if and only if $\tilde{f}=f$. Furthermore,
since the function $\theta (\ell ) = \ell \ln (\ell )$ is 
convex for $0 \leq \ell < \infty$, $D(\tilde{f},f)$ is a convex
function of $\tilde{f}$. For a fixed tolerance $c > 0$,
if ${\cal F}$ denotes the class of probability densities
over ${\mathbb R}^{n+p}$, this ensures that the "ball"
\beq
{\cal B} = \{ \tilde{f} \in {\cal F} \, : \, 
D(\tilde{f},f) \leq c \} \label{2.3}
\eeq
of densities $\tilde{f}$ located within a divergence
tolerance $c$ of the nominal density $f$ is a closed
convex set. ${\cal B}$ represents the set of all
possible true densities of random vector $z$ consistent
with the allowed mismodelling tolerance. 

Throughout this paper we shall adopt a minimax viewpoint
of robustness similar to \cite{Hub,HS}, where whenever
we seek to design an estimator minimizing an appropriately
selected loss function, a hostile player, say "nature,"
conspires to select the worst possible model in the
allowed set, here ${\cal B}$, for the performance index
to be minimized. This approach is rather conservative, 
and the performance of estimators in the presence of modelling 
uncertainties could be evaluated differently, for example by 
averaging the performance index over the entire ball 
${\cal B}$ of possible models. However this averaging
operation is computationally demanding, as it requires
a Monte Carlo simulation, and typically does not yield
analytically tractable results. It is also worth pointing 
out that the degree of conservativeness resulting from 
the selection of a minimax estimator can be controlled by 
appropriate selection of the tolerance parameter $c$,
which ensures that an adequate balance between performance
and robustness is reached. 

In this paper, we shall use the mean-square error (scaled by 
$1/2$)
\bea
J (\tilde{f},g) &=& \frac{1}{2} \tilde{E} [ || x - g(y) ||^2 ] \nnl
&=& \frac{1}{2} \int_{{\mathbb R}^{n+p}} ||x-g(y)||^2 \tilde{f} (z)dz
\label{2.4}
\eea
to evaluate the performance of an estimator $\hat{x} = g(y)$ of $x$
based on observation $y$. In (\ref{2.4}), if $v$ denotes a vector 
of ${\mathbb R}^n$ with entries $v_i$, 
\[
||v|| = (v^T v)^{1/2} = \Big( \sum_i^n v_i^2 \Big)^{1/2} 
\]
denotes the usual Euclidean vector norm. Let ${\cal G}$ denote the
class of estimators such that $\tilde{E} [ \hat{x}^2 ]$ is finite for
all $\tilde{f} \in {\cal B}$. Then the optimal robust estimator 
solves the minimax problem
\beq
\min_{g \in {\cal G}} \max_{\tilde{f} \in {\cal B}} J(\tilde{f},g) \: .
\label{2.5}
\eeq  
Since the functional $J(\tilde{f},g)$ is quadratic in $g$, and thus convex,
and linear in $\tilde{f}$, and thus concave, a saddle-point $(\tilde{f}_0,g_0)$
of minimax problem (\ref{2.5}) exists, so that
\beq
J(\tilde{f},g_0) \leq J(\tilde{f}_0,g_0) \leq J(\tilde{f}_0,g) 
\: . \label{2.6}
\eeq
However, characterizing precisely this saddle-point is difficult,
except when the nominal density is Gaussian, i.e.
\beq
f (z) \sim N( m_z , K_z ) \: , \label{2.7}
\eeq 
where in conformity with partition (\ref{2.1}) of $z$, the mean
vector $m_z$ and covariance matrix $K_z$ admit the partitions
\[
m_z = \bmat {c} 
m_x \\ m_y 
\emat \hspace*{0.15in} , \hspace*{0.15in}
K_z = \bmat {cc} 
K_x & K_{xy} \\
K_{yx} & K_y 
\emat \: . 
\]
Then it was shown in Theorem 1 of \cite{LN} (see also \cite[Sec. 7.3]{HS} for 
an equivalent result derived from a stochastic game theory perspective) that  
the saddle point of minimax poblem (\ref{2.5}) admits the following structure.
\vskip 2ex

\noindent
{\it Theorem 1:} If $f$ admits the Gaussian distribution (\ref{2.7}), the
least-favorable density $\tilde{f}_0$ is also Gaussian with distribution
\beq
\tilde{f}_0 \sim N(m_z , \tilde{K}_z) \: , \label{2.8}
\eeq
where the covariance matrix
\beq
\tilde{K}_z = \bmat {cc} 
\tilde{K}_x & K_{xy} \\
K_{yx} & K_y
\emat \label{2.9}
\eeq
is obtained by perturbing only the covariance matrix of $x$, leaving the 
cross- and co-variance matrices $K_{xy}$ and $K_y$ unchanged. Accordingly,
the robust estimator
\beq
\hat{x} = g_0 (y) = m_x + K_{xy}K_y^{-1} (y-m_y) \label{2.10} 
\eeq
coincides with the usual least-squares estimator for nominal density $f$. 
The perturbed covariance matrix $\tilde{K}_x$ can be evaluated as follows.
Let
\bea
P &=& K_x - K_{xy} K_y^{-1} K_{yx} \nnl
\tilde{P} &=& \tilde{K}_x - K_{xy} K_y^{-1} K_{yx}  \label{2.11}
\eea
denote the nominal and least-favorable error covariance matrices  
of $x$ given $y$. Then 
\beq
\tilde{P}^{-1} = P^{-1} - \lambda^{-1} I_n \: , \label{2.12}
\eeq
where $\lambda$ denotes the Lagrange multiplier corresponding to constraint
$D(\tilde{f},f) \leq c$. Note that to ensure that $\tilde{P}$ is a positive
definite matrix, we must have $\lambda > r(P)$, where $r(P)$ denotes the
spectral radius (the largest eigenvalue) of $P$.
\vskip 2ex

To explain precisely how $\lambda$ is selected to ensure that the 
Karush-Kuhn-Tucker (KKT) condition
\beq
\lambda (c - D(\tilde{f}_0,f)) = 0 \label{2.13}
\eeq
holds, observe first that for two Gaussian densities $f \sim N(m_z,K_z)$ 
and $\tilde{f} \sim N(\tilde{m}_z,\tilde{K}_z)$, the relative entropy can 
be expressed as \cite{Bas}
\beq
D(\tilde{f},f) = \frac{1}{2} \Big[ || \Delta m_z||_{K_z^{-1}}^2 +
\mbox{tr} \, \big( K_z^{-1} \tilde{K}_z - I_{n+p} \big) 
- \ln \det \big( K_z^{-1} \tilde{K}_z \big) \Big ] \: , \label{2.14}
\eeq 
where $\Delta m_z = \tilde{m}_z - m_z$ and $||v||_{K^{-1}} \defn
(v^T K^{-1} v)^{1/2}$. Then for the nominal and least-favorable
densities specified by (\ref{2.7}) and (\ref{2.8})--(\ref{2.9}),
we have $\Delta m_z =0$ and 
\bea
K_z &=& \bmat {cc}    
I_n & G_0 \\
0 & I_p
\emat \bmat {cc} 
P & 0 \\
0 & K_y 
\emat \bmat {cc}
I_n & 0 \\
G_0^T & I_p 
\emat \nnl
\tilde{K}_z &=& \bmat {cc}    
I_n & G_0 \\
0 & I_p
\emat \bmat {cc} 
\tilde{P} & 0 \\
0 & K_y 
\emat \bmat {cc}
I_n & 0 \\
G_0^T & I_p 
\emat \: , \nn
\eea
where $G_0 \defn K_{xy} K_y^{-1}$ denotes the gain matrix of the 
optimal estimator (\ref{2.10}). Then after simple algebraic
manipulations, we find 
\beq
D(\tilde{f}_0, f) = \frac{1}{2} \Big [ \mbox{tr} \,
(\tilde{P}P^{-1} -I_n ) - \ln \det (\tilde{P}P^{-1} ) \Big] 
\: . \label{2.15}
\eeq
Substituting (\ref{2.12}) gives
\beq
\gamma (\lambda) \defn D(\tilde{f}_0,f) = \frac{1}{2} \Big [ \mbox{tr}
\, ((I_n -\lambda^{-1}P)^{-1} -I_n) + \ln \det (I_n - \lambda^{-1} P) \Big ] 
\label{2.16}
\eeq 
where $\gamma (\lambda)$ is differentiable over $(r(P),\infty)$. By using
the matrix differentiation identities \cite[Chap. 8]{MN}
\bea
\frac{d~}{d\lambda} \ln \det M(\lambda) &=& \mbox{tr} \, [ M^{-1}
\frac{dM}{d\lambda} ] \nnl
\frac{d~}{d\lambda} M^{-1} (\lambda) &=& - M^{-1} \frac{dM}{d\lambda} M^{-1} \: , \nn
\eea
for a square invertible matrix function $M(\lambda)$, we find
\beq
\frac{d\gamma}{d\lambda} = -\mbox{tr} \, \big [ (I-\lambda^{-1}P)^{-1} 
\lambda^{-3}P^2 (I-\lambda^{-1}P)^{-1} \big] < 0 \label{2.17}
\eeq
so that $\gamma(\lambda)$ is monotone decreasing over $(r(P),\infty)$.
Since
\[
\lim_{\lambda \rightarrow r(P)} \gamma (\lambda ) = + \infty \hspace*{0.15in} , 
\hspace*{0.15in} \lim_{\lambda \rightarrow \infty} \gamma (\lambda ) = 0
\]
this ensures that for an arbitrary tolerance $c >0$, there exists a 
unique $\lambda >r(P)$ such that $\gamma (\lambda) = c$. 

For the case where the nominal density $f$ is non-Gaussian, some results
characterizing the solution of the minimax problem (\ref{2.5}) were
described recently in \cite{SRC}. In addition, it is worth noting
that it is assumed in Theorem 1 that the whole density $f(z) = f(x,y)$
is subject to uncertainties. But this assumption does not fit all situations. 
Consider for example a mutiple-input multiple-output (MIMO) least-squares
equalization problem for a flat channel described by the nominal 
linear model 
\beq
y = Cx + v \label{2.18}
\eeq
where $x$ denotes the transmitted data, $C$ is the channel matrix
and $v \sim N(0,R)$ represents the channel noise, which is assumed
independent of $x$. Since the transmitted data $x$ is under the control 
of the designer, its probability distribution $f(x)$ is known exactly
and it is not realistic to assume that it is affected by modelling
uncertainties. Thus if $f(y|x) \sim N(Cx,R)$ denotes the nominal
conditional distribution specified by (\ref{2.18}), the actual
density of $z$ can be represented as 
\[
\tilde{f} (z) = \tilde{f}(y|x) f(x) \: ,
\]
where $\tilde{f}(y|x)$ represents the true channel model, and where 
the data density $f(x)$ is not perturbed. This constraint changes the
structure of the minimax problem (\ref{2.5}), and a solution to this
modified problem is presented in \cite{SRC} and \cite{GL2}.

\setcounter{equation}{0}
\section{Robust Filtering Viewed as a Dynamic Game}
\label{sec:game}

We consider a robust state-space filtering problem for processes
described by a nominal Gauss-Markov state-space model of the
form
\bea
x_{t+1} &=& A_t x_t + B_t v_t \label{3.1} \\
y_t &=& C_t x_t + D_t v_t \: , \label{3.2} 
\eea
where $v_t \in {\mathbb R}^m $ is a WGN with unit variance, i.e.,
\beq
E[v_t v_s^T ] = I_m \delta(t-s) \: , \label{3.3}
\eeq
where 
\[
\delta (r) = \left \{ \begin{array} {cc}
1 & r=0 \\
0 & r \neq 0
\end{array} \right.
\]
denotes the Kronecker delta function. The noise $v_t$ is assumed to
be independent of the initial state, whose nominal distribution 
is given by 
\beq
f_0 (x_0) \sim N(m_0,P_0) \: . \label{3.4}
\eeq 
Let
\[
z_t \defn  \bmat {c}
x_{t+1} \\ y_t
\emat \: ,
\]
The model (\ref{3.1})--(\ref{3.3}) can be viewed as specifying the
nominal conditional density
\beq
\phi_t (z_t|x_t) \sim N \Big( \bmat {c}
A_t \\ C_t \emat x_t , \bmat {c}
B_t \\ D_t \emat \bmat {cc}
B_t^T & D_t^T 
\emat \Big) \: . \label{3.5}
\eeq
of $z_t$ given $x_t$. We assume that the noise $v_t$ affects all components 
of the dynamics (\ref{3.1}) and observations (\ref{3.2}), so that the covariance
matrix
\beq
K_{z_t|x_t} = 
\bmat {c}
B_t \\ D_t \emat \bmat {cc}
B_t^T & D_t^T \emat \label{3.6}
\eeq
is positive definite. To interpret this assumption, observe that in 
general, state-space models of the form (\ref{3.1})--(\ref{3.2}) 
are formed by a mixture of noisy and deterministic linear relations
(see for example the decomposition employed in \cite{LBN}). This
means that the resulting conditional densities are concentrated
on lower-dimensional subspaces of ${\mathbb R}^{n+p}$. As soon
as these subspaces are slightly perturbed, it is possible
to discriminate perfectly between the nominal and perturbed
models, i.e., the relative entropy of the two models is
infinite (the corresponding probability measures are not absolutely
continuous with respect to each other). Accordingly, when
the relative entropy is used to measure the proximity of
statistical models, all deterministic relations between dynamic
variables or observations are interpreted as immune from uncertainty,
and only relations where noise is already present can be perturbed.
Since this limitation is rather unsatisfactory, it is convenient
to assume, like earlier robust filtering studies \cite{Say,BJP,YUP},
that the noise $v_t$ affects all components of the dynamics and observations,
possibly with a very small variance for relations which are viewed as
almost certain. 

In this case, since the matrix
\[
\Gamma_t \defn \bmat{c}
B_t \\ D_t
\emat 
\]
has full row rank, we can assume without loss of generality that
$\Gamma_t$ is square and invertible, so that $m= n+p$. Otherwise,
if $m > n+p$, we can find an $m \times m$ orthonormal matrix $U_t$
which compresses the columns of $\Gamma_t$, so that
\[
\Gamma_t U_t = \bmat {cc}
\bar{\Gamma}_t & 0
\emat
\]
where $\bar{\Gamma}_t$ is invertible. Then if we denote
\[
U_t^T v_t = \bmat {c}
\bar{v}_t \\ v_t^c \emat \: ,
\]
we have
\[
\Gamma_t v_t = \bar{\Gamma}_t \bar{v}_t
\]
where $\bar{v}_t$ is a zero-mean WGN of dimension $n+p$ with unit 
covariance matrix.

Over a finite interval $0 \leq t \leq T$, the joint nominal
probability density of
\[
X_{T+1} = \bmat {c} 
x_0 \\ \vdots \\ x_t \\ \vdots \\ x_{T+1}
\emat \hspace*{0.15in} \mbox{and} \hspace*{0.15in}
Y_T = \bmat {c}
y_0 \\ \vdots \\ y_t \\ \vdots \\ y_T
\emat 
\]
can be expressed as
\beq
f(X_{T+1},Y_T) = f_0 (x_0) \prod_{t=0}^T \phi_t (z_t | x_t) \: ,
\label{3.7}
\eeq
where the initial and the combined state transition and observation densities 
are given by (\ref{3.4}) and (\ref{3.5}). Assume that the true probability
density of $X_{T+1}$ and $Y_T$ admits a similar Markov structure of the form
\beq
\tilde{f}(X_{T+1},Y_T) = \tilde{f}_0 (x_0) \prod_{t=0}^T 
\tilde{\phi}_t (z_t | x_t) \: . \label{3.8}
\eeq
By taking the expectation of
\[
\ln \Big( \frac{\tilde{f}(X_{T+1},Y_T)}{f(X_{T+1},Y_T)} \Big) = 
\ln \Big( \frac{\tilde{f}_0 (x_0)}{f_0 (x_0)} \Big) \nnl
+ \sum_{t=0}^T \Big( \frac{\tilde{\phi}_t (z_t|x_t)}{\phi_t (z_t|x_t)} \Big) 
\]
with respect to $\tilde{f}(X_{T+1},Y_T)$, we find that the relative entropy 
between $\tilde{f}(X_{T+1},Y_T)$ and $f(X_{T+1},Y_T)$ satisfies the chain
rule
\beq
D(\tilde{f},f) = D(\tilde{f}_0,f_0) + \sum_{t=0}^T D(\tilde{\phi}_t, \phi_t) \: ,
\label{3.9}
\eeq
with
\bea
D(\tilde{\phi}_t. \phi_t) &=& \tilde{E} \big[ 
\ln  \Big( \frac{\tilde{\phi}_t (z_t|x_t)}{\phi_t (z_t|x_t)} \Big) \big] \nnl 
&=& \int \int \tilde{\phi}_t (z_t|x_t) \tilde{f}_t (x_t) \ln
\Big( \frac{\tilde{\phi}_t (z_t|x_t)}{\phi_t (z_t|x_t)} \Big) dz_t dx_t \: , \label{3.10} 
\eea
where $\tilde{f}_t (x_t)$ denotes the true marginal density of $x_t$. Up to this
point, most results on robust Kalman filtering with a relative entropy
constraint have been obtained by considering a fixed interval and applying
a single constraint to the relative entropy $D(\tilde{f},f)$ of the true and  
nominal probability densities of the state and observation sequences
over the whole interval. This was also the point of view adopted in 
Section 4 of \cite{LN} which examined the robust causal Wiener filtering
problem over a finite interval. To treat the robust filtering problem over
an infinite horizon, one approach consists of dividing the divergence over 
a finite interval by the length $T$ of the interval, and letting $T$ tend
to infinity, assuming this sequence has a limit. This is the case
for stationary Gaussian processes, since in this case the limit is the
Itakura-Saito spectral distortion measure \cite[Sec. 3]{LN}. Alternatively, 
it is also possible to apply a discount factor \cite{HS1} to future additive 
terms appearing in the chain rule decomposition (\ref{3.9}). However, one 
potential weakness of applying a single divergence constraint to the
filtering problem over a finite or infinite interval is that it allows
the maximizer (nature) to identify the moment where the dynamic model
(\ref{3.1})-(\ref{3.3}) is most susceptible to distortions and to
allocate most of the distortion budget specified by the tolerance 
parameter $c$ to this single element of the model. If the purpose
of robust filtering is to protect the estimator from modelling
inaccuracies, and if the modeller exercises the same level of effort 
to characterize each time component of the model (\ref{3.1})--(\ref{3.3}), 
it may be more appropriate to specify separate modelling tolerances
for each time step of the transition density (\ref{3.5}). Such a 
viewpoint has in fact been adopted widely \cite{PS,GC,Say} in the
robust state-space filtering literature, except that in these earlier 
studies the tolerance is usually expressed in terms of matrix bounds 
involving the matrices $A_t$, $B_t$, $C_t$ and $D_t$ parametrizing the 
state dynamics and observations. The main difference with these earlier 
studies is that we use here the relative entropy $D(\tilde{\phi}_t,\phi_t)$ 
between the true and nominal transition and observation densities 
$\tilde{\phi}_t (z_t|x_t)$ and $\phi_t (z_t|x_t)$ at time $t$ to 
measure modelling errors. 

The expression (\ref{3.10}) for the relative entropy raises immediately
the issue of how to choose the  probability density $\tilde{f}_t (x_t)$
used to evaluate the divergence. We assume that, like the estimating
player, at time $t$ the maximizer has access to the observations $\{ 
y(s), ~ 0 \leq s \leq t-1 \}$ collected prior to this point. In addition,
it is reasonable to hold the maximizer to the same Markov structure
(specified by (\ref{3.8})) as the estimating player. Therefore, the
maximizer is required to commit to all the least-favorable model 
components $\tilde{\phi}_s (z_s | x_s)$ with $0 \leq s \leq t-1$ 
generated at earlier stages of its minimax game with the estimating 
player. Using the terminology coined in \cite{HS1,HS}, the maximizer 
operates "under commitment." Thus if $Y_{t-1}$ denotes the vector 
formed by the observations $\{ y_s ,~ 0 \leq s \leq t-1 \}$, we use
the conditional density $\tilde{f}_t (x_t|Y_{t-1})$ based on the
least favorable model and the given observations prior to time $t$,
to evaluate the divergence (\ref{3.10}) between the true and nominal transition
and observation densities. The model mismatch tolerance can therefore
be expressed as
\beq
\tilde{E} \Big [ \ln \Big( \frac{\tilde{\phi}_t (z_t|x_t)}{\phi_t (z_t|x_t)} \Big)
\Big | Y_{t-1} \Big ] \leq c_t \: , \label{3.11}
\eeq
where $c_t$ denotes the tolerance parameter for the time $t$ component of the
model, with
\bea
\lefteqn{ \tilde{E} \Big[ \ln \Big ( \frac{\tilde{\phi}_t (z_t|x_t)}{\phi_t (z_t|x_t)}
\Big) \Big | Y_{t-1} \Big ]} \nnl
&=& \int \int \tilde{\phi}_t (z_t|x_t) \tilde{f}_t (x_t|Y_{t-1})
\ln \Big( \frac{\tilde{\phi}_t (z_t|x_t)}{\phi_t (z_t|x_t)} \Big) dz_t dx_t \: . 
\label{3.12}
\eea

Let ${\cal B}_t$ denote the convex ball of functions $\tilde{\phi}_t (z_t|x_t)$
satisfying inequality (\ref{3.11}). If ${\cal G}_t$ denotes the class of
estimators with finite second-order moments with respect to all densities
$\tilde{\phi}_t (z_t|x_t)\tilde{f}_t (x_t|Y_{t-1})$ such that $\tilde{\phi}_t (z_t|x_t)
\in {\cal B}_t$, the dynamic minimax game we consider can be expressed as
\beq
\min_{g_t \in {\cal G}_t} \max_{\tilde{\phi}_t \in {\cal B}_t}
J_t (\tilde{\phi}_t,g_t) \label{3.13}
\eeq
where
\bea
J_t (\tilde{\phi}_t,g_t) &=& \frac{1}{2} \tilde{E} \Big[ ||x_{t+1} - g_t (y_t)||^2|
Y_{t-1} \Big] \nnl
&=& \frac{1}{2} \int \int || x_{t+1} - g_t (y_t)||^2 \tilde{\phi}_t (z_t|x_t) \tilde{f}_t (x_t|
Y_{t-1}) dz_t dx_t \: . \label{3.14}
\eea
denotes the mean-square error of estimate $\hat{x}_{t+1} = g_t (y_t)$ of $x_{t+1}$
evaluated with respect to the true probability density of $z_t$. Note that
since $\hat{x}_{t+1}$ is a function of $Y_t$, it depends not only on $y_t$,
but also on earlier observations, but this dependency is suppressed
to simplify notations. Note that by taking iterated expectations, we have
\bea
\tilde{E} [ || x_{t+1} - \hat{x}_{t+1} ||^2] &=& \tilde{E} [ \tilde{E} [ ||x_{t+1} -
\hat{x}_{t+1} ||^2 \ Y_{t-1} ] ] \nnl
&=& \int J_t ( \tilde{\phi}_t , g_t ) \tilde{f}_{t-1} (Y_{t-1}) d Y_{t-1} \: , \label{3.15}
\eea
where $\tilde{f}_{t-1} (Y_{t-1})$ represents the least favorable density of
observation vector $Y_{t-1}$ based on the least-favorable model increments 
selected by the maximizer up to time $t-1$. Since this density is non-negative
and independent of both $\tilde{\phi}_t$ and $g_t$, the game (\ref{3.13}) 
is equivalent to 
\beq
\min_{g_t \in {\cal G}_t} \max_{\tilde{\phi}_t \in {\cal B}_t}
\tilde{E} [ || x_{t+1} - \hat{x}_{t+1} ||^2 ] \: . \label{3.16}
\eeq
This indicates that the robust least-squares filtering problem we consider
is local in the sense that the estimator and the maximizer focus
respectively on minimizing and maximizing the mean-square estimation
error at the current time. In other words, in addition to being committed
to past least-favorable model increments $\tilde{\phi}_s$ with $0 \leq s 
\leq t-1$ it has already selected, the maximizer is confined to a {\it myopic
strategy}, where $\tilde{\phi}_t$ is selected exclusively to maximize the
mean-square estimation error at time $t$. By doing so, the maximizer 
foregoes the possibility of trading off a lesser increase in the mean-square
error at the current time against larger increases in the future.  

If we compare the dynamic estimation game (\ref{3.13}) and its static
counterpart (\ref{2.5}), we see that the two problems are similar, but
the dynamic game (\ref{3.13}) includes a conditioning operation with
respect to the prior state $x_t$, combined with an averaging operation
with respect to $\tilde{f}_t (x_t |Y_{t-1})$. Thus Lemma 1 and 
Theorem 1 of \cite{LN} need to be extended slightly to accommodate 
these differences. Before proceeding with this task it is worth pointing
out that we do not require that $\tilde{\phi}_t (z_t|x_t)$ should
be a conditional probability density for each $x_t$. It is only required
that the product $\tilde{\phi}_t (z_t|x_t) \tilde{f}_t (x_t|Y_{t-1})$
should be a probability density for
\[
\bmat {c}
z_t \\ x_t 
\emat = \bmat {c}   
x_{t+1} \\ y_t \\ x_t 
\emat \: , 
\]
so that
\beq
I_t (\tilde{\phi}_t) \defn \int \int \tilde{\phi}_t (z_t|x_t) \tilde{f}_t (x_t|Y_{t-1}) dz_t dx_t
= 1 \: . \label{3.17}
\eeq
To put it another way, the maximizer's commitment to earlier components
of the least-favorable model is only of an a-priori nature, since the
a-posteriori marginal density of $x_t$ specified by the joint density
$\tilde{\phi}_t (z_t|x_t) \tilde{f}_t (x_t|Y_{t-1})$ is not
required to coincide with the a priori density $\tilde{f}_t (x_t|Y_{t-1})$.  
However, by integrating out $z_t$, the resulting a-posteriori density
will be of the form $\psi_t (x_t) \tilde{f}_t (x_t |Y_{t-1})$ where
$\psi_t (x_t)$ does not depend on $Y_{t-1}$.  
\vskip 2ex

\noindent
{\bf Relation to prior work:} At this point it is possible to compare
precisely the robust filtering problem discussed here with earlier 
formulations of robust filtering with a relative entropy tolerance 
considered in \cite{HS,HS1,LN,YUP}. In \cite{HS,HS1}, Hansen and Sargent 
introduce the likelihood ratio function
\beq
M_t = \frac{\tilde{f} (X_{t+1},Y_t)}{f(X_{t+1},Y_t)} \label{3.18}   
\eeq
between the distorted joint density of states $\{ x_s, 0 \leq s \leq t+1\}$ 
and observations $\{ y_s, 0 \leq s \leq t \}$ and their nominal joint
density. Setting
\[
M_{-1} = \frac{\tilde{f}_0 (x_0)}{f_0 (x_0)} \: ,
\]
if ${\cal Z}_t$ denotes the $\sigma$-field generated by $\{ z_s, 
0 \leq s \leq t \}$ and $x_0$, we have
\[
E[ M_t | {\cal Z}_{t-1} ] = M_{t-1} \: ,
\]
for $t \geq 0$, where $E[ \cdot ]$ denotes the expectation with respect to the 
nominal model distribution, so $M_t$ is a martingale. Then the
relative entropy between the perturbed and nominal models
over interval $0 \leq t \leq T$ can be expressed as 
\beq
D(\tilde{f},f) = E[ M_T \ln M_T ] \: . \label{3.19} 
\eeq
Let
\beq
J = \sum_{t=0}^T \tilde{E} [||x_t - \hat{x}_t ||^2 ] 
= \sum_{t=0}^T E[ M_t ||x_t - \hat{x}_t||^2 ]
\label{3.20}
\eeq
denote the sum of the mean square estimation errors over
interval $[0,T]$, where the estimate $\hat{x}_t = g_t (Y_t)$ depends
causally on the past and current observations. If $g = (g_t, t \geq 0)$, 
the robust Kalman filtering problem of \cite{HS,HS1} and the
robust causal Wiener filtering problem of \cite[Sec. 4]{LN} 
can be written as 
\beq
\min_{g} \max_{\tilde{f}} J \label{3.21}
\eeq
where $\tilde{f}$ satisfies the constraint
\beq
D(\tilde{f},f) = E[M_T \ln M_T ] \leq c \: . \label{3.22}
\eeq
Since the causal Wiener filtering problem of \cite[Sec. 4]{LN}
is treated as a constrained static estimation problem, it does not 
require any additional structure. On the other hand for the robust Kalman
filtering problem of \cite{HS,HS1}, the minimax problem (\ref{3.21})
is formulated dynamically by introducing the martingale increments
\beq
m_t = \left \{ \begin{array} {ccc}
\frac{M_t}{M_{t-1}} & \mathrm{if} &  M_{t-1} > 0 \\[1ex]
1 & \mathrm{if} &  M_{t-1} = 0 \: ,
\end{array} \right. \label{3.23}
\eeq
which can be expressed here as
\[
m_t = \frac{\tilde{\phi}(z_t|x_t)}{\phi(z_t|x_t)} \: .
\]
Then the chain rule (\ref{3.10}) takes the form
\beq
E[M_T \ln M_T ] = \sum_{t=0}^T E[ M_{t-1} E[ m_t \ln m_t |{\cal Z}_{t-1}] ] + 
E[ M_{-1} \ln M_{-1}] \: , \label{3.24}
\eeq
and the dynamic minimax problem considered in \cite{HS,HS1} can be expressed
as 
\beq
\min_{g_t , t \geq 1} \, \max_{ m_t, t \geq 1 } J  \label{3.25}
\eeq
where $g_t$ and $m_t \geq 0$ are adapted respectively to ${\cal Y}_t$ 
(the sigma field spanned by $Y_t$) and ${\cal Z}_t$. Comparing 
(\ref{3.16}) and (\ref{3.25}), we see that (\ref{3.16}) represents
just a local or incremental version of minimax problem (\ref{3.25}).
In this respect it is worth pointing out that by integrating
the incremental relative entropy constraint (\ref{3.11}) with
respect to the least favorable density $\tilde{f}_{t-1} (Y_{t-1})$
of $Y_{t-1}$ and summing over all $0 \leq t \leq T$, we find that
cumulatively, the incremental constraints (\ref{3.11}) imply
\beq
D(\tilde{f},f) \leq  \sum_{t=0}^T c_t + D(\tilde{f}_0, f_0) \: , 
\label{3.26} 
\eeq
which has the form (\ref{3.22}). In other words, the incremental
constraint (\ref{3.11}) just represents a way of dividing
the relative entropy tolerance budget $c$ in separate 
portions $c_t$ allocated to the distortion of the Markov
model transition at each time step.

The robust filtering problem we consider can also be interpreted in
terms of the robust filtering approach described in \cite{Say}. To do so, 
assume that the distorted transition 
\[
\tilde{\phi}_t (z_t|x_t) \sim N \Big( \bmat {c}
\tilde{A}_t \\ \tilde{C}_t 
\emat x_t , \tilde{K}_{z_t|x_t} \Big)
\]
is Gaussian and admits a parametrization similar to (\ref{3.5}). 
By using expression (\ref{2.14}) for the relative entropy 
of Gaussian densities, and denoting 
\[
\Delta A_t \defn  \tilde{A}_t - A_t  \hspace*{0.1in} , \hspace*{0.1in} 
\Delta C_t \defn  \tilde{C}_t - C_t : ,   
\]
we find that conditioned on the knowledge of $x_t$ 
\beq
D(\tilde{\phi}_t (.|x_t) , \phi_t (.|x_t) ) = \frac{1}{2} 
\Big[ || \bmat{c} 
\Delta A_t \\
\Delta C_t 
\emat x_t ||_{K_t^{-1}}^2 + \mbox{tr} \Big( K_t^{-1} \tilde{K}_t -I \Big) 
- \ln \det ( K_t^{-1} \tilde{K}_t ) \Big ] \: , \label{3.27}
\eeq
where for simplicity we fave used the compact notation $K_t = K_{z_t|x_t}$
and $\tilde{K}_t = \tilde{K}_{z_t|x_t}$. Then, if
\[
\tilde{f}_t (x_t |Y_{t-1}) \sim N(\hat{x}_t , V_t) \: ,
\]
is Gaussian and $W_t \defn V_t + \hat{x}_t \hat{x}_t^T$, the expression (\ref{3.10}) 
for the relative entropy between the distorted and nominal transitions at 
time $t$ yields
\beq
D(\tilde{\phi}_t,\phi_t ) = \frac{1}{2} 
\Big[ || K_t^{-1/2} \bmat{c} 
\Delta A_t \\
\Delta C_t 
\emat W_t^{1/2} ||_F^2 + \mbox{tr} \Big( K_t^{-1} \tilde{K}_t -I \Big) 
- \ln \det ( K_t^{-1} \tilde{K}_t ) \Big ] \: , \label{3.28}
\eeq
where $K_t^{1/2}$ and $W_t^{1/2}$ are arbitrary matrix square roots of $K_t$ 
and $W_t$, respectively,  and $||.||_F$ denotes the Frobenius norm of a 
matrix \cite[p. 291]{HJ}. The first term of expression (\ref{3.28}) is 
a weighted matrix norm of the perturbations $\Delta A_t$ and $\Delta C_t$ 
of the state transition dynamics which is similar in nature to the
mismodelling measures considered in \cite{Say,GC}. On the other hand,
the second term models the distortion of the process and measurement
noise covariance $K_t = K_{z_t|x_t}$, and is different from the 
distortion metrics considered in \cite{Say,GC}. Thus the robust
filtering problem we consider can on one hand be viewed as
incremental version of the results of \cite{HS,HS1} for robust 
filtering with a relative entropy tolerance, but it can also 
be viewed as a variant of the robust Kalman filters discussed
in \cite{Say,GC} with a different local mismodelling measure. 
\vskip 2ex

\noindent
{\bf Formulation without commitment:} The commitment assumption
can be removed from the incremental formulation of robust filtering
described above by adopting the conceptual framework of Hansen and
Sargent in \cite[Chap. 18]{HS} and \cite{HS2,HS3} for robust filtering
without commitment. The key idea is, at time $t$, to apply 
distortions to both the transition dynamics $\phi _t (z_t|x_t)$ 
and the least-favorable a-priori distribution $\tilde{f}_t (x_t|Y_{t-1})$.
In \cite{HS2} this is accomplished by introducing two distortion
operators with constant risk-sensitivity parameters $\theta_1$  
and $\theta_2$, acting respectively on the transition dynamics 
and on the a-priori density at time $t$ based on the prior
distortions and observations. Here, if $\check{f}_t (x_t |Y_{t-1})$
denotes a distorted version of $\tilde{f}_t (x_t|Y_{t-1})$, in
addition to the constraint (\ref{3.11}) for transition dynamics
distortion, we could impose a constraint of the form
\beq
D(\check{f}_t,\tilde{f}_t ) \leq d_t \label{3.29}
\eeq
on the allowed distortion of the least favorable conditional 
distortion for $x_t$ based on the past observations $Y_{t-1}$
Then in the minimax problem (\ref{3.13}) the maximization can
be performed jointly over pairs $(\tilde{\phi}_t, \check{f}_t)$
satisfying constraints (\ref{3.11}) and (\ref{3.29}). Since
the two constraints are convex, the structure of the resulting
minimax problem is similar to (\ref{3.13}). The main difficulty
is algorithmic. The inner maximization introduces two Lagrange
multipliers which need to be selected such that 
Karush-Kuhn-Tucker conditions are satisfied. The computation of
the Lagrange multipliers seems rather difficult, in contrast to the 
case of a single Lagrange multiplier arising from the formulation
of robust filtering with commitment. So we focus here on the
case with commitment, leaving open the possibility that an 
implementable algorithm might be developed later for robust
filtering without commitment under incremental constraints
(\ref{3.11}) and (\ref{3.29}). Finally, note that a third
option would be to combine the distortions for the
transition density $\phi_t (z_t|x_t)$ and for conditional
density $\tilde{f}_t (x_t|Y_{t-1})$ and to apply a
single relative entropy constraint to the product distorted
density $\tilde{\phi}_t (z_t|x_t) \check{f}_t (x_t|Y_{t-1})$.
Since this product corresponds to the least favorable joint
demsity of $z_t$ and $x_t$, Theorem 1 is applicable to
this problem, so it is not necessary to analyze this version 
of the robust filtering problem.
\vskip 2ex

\setcounter{equation}{0}
\section{Robust Minimax Filter}
\label{sec:filter}

The solution of the dynamic game (\ref{3.13}) relies on extending Lemma 1
of \cite{LN} to the dynamic case. We start by observing that the objective
function $J_t (\tilde{\phi}_t, g_t)$ specified by (\ref{3.14}) is 
quadratic in $g_t$, and thus convex, and linear in $\tilde{\phi}_t$ and
thus concave. The set ${\cal B}_t$ is convex and compact. Similarly ${\cal G}_t$
is convex. It can also be made compact by requiring that the second moment
of estimators $g_t \in {\cal G}_t$ should have a fixed but large upper bound. 
Then by Von Neumann's minimax theorem \cite[p. 319]{AE}, there exists a saddle
point $(\tilde{\phi}_t^0 , g_t^0)$ such that
\beq
J_t (\tilde{\phi}_t, g_t^0) \leq J_t (\tilde{\phi}_t^0,g_t^0) \leq
J_t (\tilde{\phi}_t^0, g_t) \: . \label{4.1} 
\eeq 
The real challenge is, however, not to establish the existence of a saddle
point, but to characterize it completely. The second inequality in (\ref{4.1})
implies that estimator $g_t^0$ is the conditional mean of $x_{t+1}$ given
$Y_t$ based on the least-favorable density
\beq
\tilde{f}_{t+1} (x_{t+1}|Y_t) = \frac{\int \tilde{\phi}_t^0 (z_t|x_t)
\tilde{f}_t (x_t|Y_{t-1}) dx_t}{\int \int \tilde{\phi}_t^0 (z_t|x_t)
\tilde{f}_t (x_t|Y_{t-1}) dx_{t+1}dx_t} \label{4.2}
\eeq
obtained by marginalization and application of Bayes' rule to the least-favorable
joint density $\tilde{\phi}_t^0 (z_t|x_t) \tilde{f}_t (x_t|Y_{t-1})$
of $(z_t,x_t) = (x_{t+1},y_t,x_t)$ given $Y_{t-1}$. The robust estimator
is then given by
\bea
\hat{x}_{t+1} = g_t^0 (y_t) &=& \tilde{E} [ x_{t+1} |Y_t ] \nnl
&=&  \int x_{t+1} \tilde{f}_{t+1} (x_{t+1}|Y_t) dx_{t+1} \: . \label{4.3}
\eea
Together, the conditional density evaluation (\ref{4.2}) and expectation
(\ref{4.3}) implement the second inequality of saddle point identity (\ref{4.1}).
Let us turn now to the first inequality. For a fixed estimator $g_t^0$, it requires
finding the joint transition and observation density $\tilde{\phi}_t^0 (z_t|x_t)$
maximizing $J_t(\tilde{\phi}_t,g_t^0)$ under the divergence constraint (\ref{3.11}).
The solution of this problem takes the following form.
\vskip 2ex

\noindent
{\it Lemma 1:} For a fixed estimator $g_t \in {\cal G}_t$, the function 
$\tilde{\phi}_t^0$ maximizing $J_t (\tilde{\phi}_t,g_t)$ under constraints   
(\ref{3.11}) and (\ref{3.17}) is given by
\beq
\tilde{\phi}_t^0 (z_t|x_t) = \frac{1}{M_t (\lambda_t)} \exp \Big (\frac{1}{
2\lambda_t} ||x_{t+1} -g_t(y_t)||^2 \Big) \phi_t (z_t | x_t) \: . \label{4.4}
\eeq
In this expression, the normalizing constant $M_t (\lambda_t )$ is selected such
that (\ref{3.17}) holds. Furthermore, given a tolerance $c_t > 0$, there
exits a unique Lagrange multiplier $\lambda_t >0$ such that
\beq
D_t (\tilde{\phi}_t^0,\phi_t ) = c_t \: . \label{4.5}
\eeq
\vskip 2ex

\noindent
{\it Proof:} For a given $g_t$, the function $J_t (\tilde{\phi}_t,g_t)$
is linear in $\tilde{\phi}_t$ and thus concave over the closed convex set 
${\cal B}_t$, so it admits a unique maximum in ${\cal B}_t$. Because of the
linearity of $J_t$ with respect to $\tilde{\phi}_t$, this maximum is
in fact located on the bpundary of ${\cal B}_t$. To find the maximum, consider 
the Lagrangian
\beq
L_t (\tilde{\phi}_t,\lambda_t ,\mu_t ) = J_t(\tilde{\phi}_t,g) + \lambda_t
(c_t - D_t(\tilde{\phi}_t,\phi_t)) + \mu_t ( 1 -I_t (\tilde{\phi}_t )) \: ,
\label{4.6}
\eeq
where the Lagrange multipliers $\lambda_t \geq 0$ and $\mu_t$ are associated to 
inequality constraint (\ref{3.11}) and equality constraint (\ref{3.17}),
respectively. We do not require explicitly that $\tilde{\phi}_t (z_t|x_t)$
should be nonnegative, since the form (\ref{4.4}) of the maximizing solution
indicates that this constraint is satisfied automatically.

Then the Gateaux derivative \cite[p. 17]{BNO} of $L_t$ with respect to 
$\tilde{\phi}_t$ in the  direction of an arbitrary function $u$ is given by
\bea
\lefteqn{\nabla_{\tilde{\phi}_t, u} L_t (\tilde{\phi}_t,\lambda_t ,\mu_t ) =
\lim_{h \rightarrow 0} \frac{1}{h} [ L_t (\tilde{\phi}_t + hu,\lambda_t ,\mu_t ) 
- L_t (\tilde{\phi}_t,\lambda_t ,\mu_t ) ]} \nnl
&=& \int \int \Big [ \frac{1}{2} ||x_{t+1} -g_t (y_t)||^2 - (\lambda_t + \mu_t) -
\lambda_t \ln \big( \frac{\tilde{\phi}_t}{\phi_t} \big) \Big] u(z_t,x_t)  \tilde{f}_t (x_t|Y_{t-1})
dz_t dx_t \: . \hspace*{0.2in} \label{4.7}
\eea  
The Lagrangian is maximized by setting $\nabla_{\tilde{\phi}_t , u} L_t (\tilde{\phi}_t,
\lambda_t,\mu_t) = 0$ for all functions $u$. Assuming $\lambda_t >0$, this gives
\beq
\ln \Big( \frac{\tilde{\phi}_t}{\phi_t} \Big) = \frac{1}{2\lambda_t}   
|| x_{t+1} - g_t (y_t)||^2 - \ln M_t \: , \label{4.8}
\eeq
where
\[
\ln M_t \defn 1 + \frac{\mu_t}{\lambda_t} \: .
\]
Exponentiating (\ref{4.8}) gives (\ref{4.4}), where to ensure that normalization
(\ref{3.14}) holds, we must select
\beq
M_t (\lambda_t) = \int \int \exp \Big( \frac{1}{2\lambda_t} ||x_{t+1} - g_t (y_t)||^2 
\Big) \phi_t (z_t|x_t) \tilde{f}_t (x_t|Y_{t-1}) dz_t dx_t \: . \label{4.9}
\eeq
At this point, all what is left is finding a Lagrange multiplier $\lambda_t$ such that
the solution $\tilde{\phi}_t^0$ given by (\ref{4.4}) satisfies the KKT condition 
\[
\lambda_t (c_t - D_t (\tilde{\phi}_t^0, \phi_t) = 0 \: .
\]
Since we already know that the maximizing $\tilde{\phi}_t^0$ is on the boundary
of ${\cal B}_t$, the Lagrange multiplier $\lambda_t >0$, so the KKT condition
reduces to (\ref{4.5}). By substituting (\ref{4.4}) inside expression 
(\ref{3.11}) for $D(\tilde{\phi}_t^0,\phi_t)$, we find
\beq
D(\tilde{\phi}_t^0,\phi_t) = \frac{1}{\lambda_t} J_t (\tilde{\phi}_t^0,g_t)  
-\ln (M_t (\lambda_t)) \: . \label{4.10}  
\eeq
Differentiating $\ln M_t (\lambda_t)$ gives
\beq
\frac{d~}{d\lambda_t}\ln M_t (\lambda_t)  = -\frac{1}{\lambda_t^2} J_t (\tilde{\phi}_t^0,g_t) \: ,
\label{4.11} 
\eeq
so that
\beq
\gamma (\lambda_t) \defn D(\tilde{\phi}_t^0,phi_t) = -\lambda_t \frac{d~}{d\lambda_t} 
\ln M_t (\lambda_t) -\ln (M_t (\lambda_t) \: . \label{4.12} 
\eeq 
The derivative of $\gamma (\lambda_t )$ is given by
\bea
\frac{d\gamma_t}{d\lambda_t} &=& -\lambda_t \frac{d^2~}{d\lambda_t^2} \ln M_t 
-2 \frac{d~}{d\lambda_t} \ln M_t \nnl
&=& -\frac{1}{\lambda_t} \Big[ \frac{d~}{d\lambda_t} \Big( \lambda_t^2 
\frac{d~}{d\lambda_t} \ln M_t \Big) \Big] = \frac{1}{\lambda_t} \frac{d~}{d\lambda_t} 
J_t (\tilde{\phi}_t^0, g_t) \nnl
&=& -\frac{1}{4\lambda_t^3} \tilde{E} \big [ \big( ||x_{t+1} - g_t (y_t)||^2 
- \tilde{E} [ ||x_{t+1} -g_t (y_t) ||^2 | Y_{t-1}] \big)^2 | Y_{t-1} \big] < 0 \: , \label{4.13} 
\eea
so that $\gamma (\lambda_t)$ is a monotone decreasing function of $\lambda_t$. As
$\lambda_t \rightarrow \infty$ we have obviously $\tilde{\phi}_t^0 \rightarrow \phi_t$,
so that $\gamma (\infty) =0$. Thus provided $c_t$ is located in the range of
$\gamma_t$, which is the case if $c_t$ is sufficiently small, there exists
a unique $c_t$ such that $\gamma_t (\lambda_t ) = c_t$. \hfill $\Box$ 
\vskip 2ex

Note that Lemma 1 makes no assumption about the form of the nominal transition
density $\phi_t (z_t|x_t)$ and a priori density $\tilde{f}_t (x_t|Y_{t-1})$.
Without additional assumptions, it is difficult to characterize precisely
the range of function $\gamma_t$. When both of these densities are Gaussian,
it will be shown below that the range of $\gamma_t$ is ${\mathbb R}^+$,
so that any positive divergence tolerance $c_t$ can be achieved. However,
in practice the tolerance $c_t$ needs to be rather small in order to
ensure that the robust estimator is not overly conservative. At this
point is is also worth observing that Lemma 1 is just a variation 
of Theorem 2.1 in \cite[p. 38]{Kul} which sought to construct the minimum 
discrimination density (i.e., the density minimizing the divergence)
with respect to a nominal density under various moment constraints. Here we
seek to maximize the moment $\tilde{E}[ ||x_{t+1} - g_t (y_t)||^2 |Y_{t-1}]$ 
under a divergence constraint. From an optimization point of view, the two 
problems are obviously similar, and in fact the functional form (\ref{4.4}) 
of the solution is the same for both problems.
 
Up to this point we have made no assumption on either the nominal transition
and observations density $\phi_t (z_t|x_t)$ and estimator $g_t$, and in
the characterization of the saddle point solution $(\tilde{\phi}_t^0 ,g_t^0 )$
provided by identities (\ref{4.3}) and (\ref{4.4}), the robust estimator
$g_t^0$ depends on least-favorable transition function $\tilde{\phi}_t^0$,
and the least favorable transition density $\tilde{\phi}_t^0$ depends 
on robust estimator $g_t^0$. This type of deadlock is typical of
saddle point analyses, and to break it, we introduce now the assumption
that $\phi_t (z_t|x_t)$ admits the Gaussian form (\ref{3.5}) where as indicated
earlier, the covariance matrix $K_{z_t|x_t}$ is positive definite, and we assume 
also that at time $t$ the a-priori conditional density
\beq
\tilde{f}_t (x_t|Y_{t-1}) \sim {\cal N} (\hat{x}_t,V_t) \: . \label{4.14} 
\eeq
Then, observe that the distortion term $\exp( ||x_{t+1}-g_t(y_t)||^2/(2\lambda_t))$
appearing in expression (\ref{4.4}) for the least-favorable transition function
$\tilde{\phi}_t^0 (z_t|x_t)$ depends only on $z_t$, but not $x_t$. Accordingly,
if we introduce the marginal densities
\bea
\bar{f}_t (z_t|Y_{t-1}) &=& \int \phi_t (z_t|x_t) \tilde{f}_t (x_t |Y_{t-1})  
dx_t \label{4.15} \\
\tilde{f}_t^0 (z_t|Y_{t-1}) &=& \int \tilde{\phi}_t^0 (z_t|x_t)
\tilde{f}_t (x_t|Y_{t-1}) dx_t \: , \label{4.16}
\eea
the density $\bar{f} (z_t|Y_{t-1})$ can be viewed as the pseudo-nominal
density of $z_t = (x_{t+1},y_t)$ conditioned on $Y_{t-1}$ computed from 
the conditional least favorable density $\tilde{f}_t (x_t|Y_{t-1})$ and 
nominal transition density $\phi (z_t|x_t)$, and from (\ref{4.4}) we obtain
\beq
\tilde{f}_t^0 (z_t|Y_{t-1}) = \frac{1}{M(\lambda_t)} \exp \Big(
\frac{1}{2\lambda_t} ||x_{t+1} - g_t (y_t)||^2 \Big) \bar{f}_t (z_t|Y_{t-1})
\: . \label{4.17} 
\eeq
Since densities $\phi_t (z_t|x_t)$ and $\tilde{f}_t (x_t|Y_{t-1})$ are
both Gaussian, the integration (\ref{4.15}) yields a Gaussian pseudo-nominal 
density
\beq
\bar{f}_t (z_t|Y_{t-1}) \sim {\cal N} \Big( \bmat {c} 
A_t \\ C_t \emat \hat{x}_t , K_{z_t} ) \label{4.18}
\eeq
where the conditional covariance matrix $K_{z_t}$ is given by 
\beq
K_{z_t} = \bmat {c} 
A_t \\ C_t \emat V_t \bmat {cc}
A_t^T & C_t^T \emat + \bmat {c}
B_t \\ D_t \emat \bmat {cc}
B_t^T & D_t^T \emat \: . \label{4.19}  
\eeq
By integrating out $x_t$ in (\ref{4.9}), we find 
\[
M_t (\lambda_t ) = \int \exp \Big( \frac{1}{2\lambda_t} ||x_{t+1} -g_0(y_t)||^2 \Big)
\bar{f}_t (z_t|Y_{t-1}) dz_t \: , 
\]
which ensures that $\tilde{f}_t^0 (z_t|Y_{t-1})$ is a probability density. 
Furthermore, by direct substitution, we have
\beq
D(\tilde{f}_t^0,\bar{f}_t) = D(\tilde{\phi}_t^0,\phi_t) = c_t \: . \label{4.20}
\eeq 
The least-favorable density $\tilde{f}_{t+1} (x_{t+1}|Y_t)$ specified by 
(\ref{4.2}) can also be expressed in terms of $\tilde{f}_t^0 (z_t|Y_{t-1})$
as
\beq
\tilde{f}_{t+1} (x_{t+1}|Y_t) = \frac{\tilde{f}_t^0 (z_t|Y_{t-1})}{
\int \tilde{f}_t^0 (z_t|Y_{t-1}) dx_{t+1}} \: . \label{4.21}
\eeq
\vskip 2ex

\noindent
{\bf Equivalent Static Problem:} Let 
\beq
\bar{\cal B}_t = \{ \tilde{f}_t : D(\tilde{f}_t,\bar{f}_t) \leq c_t \} \label{4.22}
\eeq
denote the ball of distorted densities $\tilde{f}_t (z_t)$ within a divergence
tolerance $c_t$ of pseudo-nominal density $\bar{f}_t (z_t|Y_{t-1})$.
To this ball we can of course attach a static minimax estimation problem
\beq
\min_{g_t \in {\cal G}_t} \max_{\tilde{f}_t \in \bar{\cal B}_t} J_t (\tilde{f}_t,g_t) 
\: . \label{4.23} 
\eeq
The solution of this problem satisfes the saddle point inequality
\beq
J( \tilde{f}_t, g_t^0) \leq J(\tilde{f}_t^0,g_t^0) \leq J_t (\tilde{f}_t^0,g_t) \: .
\label{4.24}
\eeq
At this point, observe that if $(\tilde{\phi}_t^0,g_t^0)$ solves
the dynamic minimax game (\ref{3.12}), and if $\tilde{f}_t^0$ is given
by (\ref{4.17}) with $g_t=g_t^0$, where $\lambda_t$ is selected such that 
constraint (\ref{4.20}) is satisfied, then $(\tilde{f}_t^0,g_t^0)$ is a saddle 
point of the static problem (\ref{4.23}). In other words, the marginalization
operation (\ref{4.16}) has the effect of mapping the solution of dynamic 
game (\ref{3.12}) into a solution of the static estimation problem 
problem (\ref{4.23}). Note indeed that the solution of the maximization
problem formed by the first inequality of (\ref{4.24}) is given by (\ref{4.17})
with $g_t = g_t^0$. Similarly, since $g_t^0$ is the mean of the conditional
density $\tilde{f}_{t+1} (x_{t+1}|Y_t)$ specified by (\ref{4.21}), it 
obeys the second inequality of (\ref{4.24}).
   
Since the pseudo-nominal density $\bar{f}_t (z_t|Y_{t-1})$ specifying the
center of ball $\bar{\cal B}_t$ is Gaussian, Theorem 1 is applicable with
$f \rightarrow \bar{f}_t$, $\tilde{f}_0 \rightarrow \tilde{f}_t^0$ and $g_0 
\rightarrow g_t^0$. Hence the least-favorable density takes the form
\beq
\tilde{f}_t^0 (z_t |Y_{t-1}) \sim {\cal N} \Big( \bmat {c}
A_t \\ C_t \emat \hat{x}_t , \tilde{K}_{z_t} \Big ) 
\label{4.25}
\eeq
where the covariance matrix
\[
\tilde{K}_{z_t} = \bmat {cc} 
\tilde{K}_{x_{t+1}} & K_{x_{t+1} y_t} \\
K_{y_t x_{t+1}} & K_{y_t}
\emat
\]
is obtained by perturbing only the (1,1) block 
\[
K_{x_{t+1}} = A_t V_t A_t^T + B_t B_t^T 
\]
of the covariance matrix $K_{z_t}$ given by (\ref{4.19}). The robust estimator
takes the form 
\beq
\hat{x}_{t+1} = g_t^0 (y_t) = A_t \hat{x}_t + G_t (y_t - C_t \hat{x}_t)
\label{4.26}
\eeq
with the matrix gain
\beq 
G_t = K_{x_{t+1}y_t}K_{y_t}^{-1} = (A_t V_t C_t^T + B_t D_t^T)(C_t V_t C_t^T 
+ D_t D_t^T)^{-1} \: . \label{4.27}
\eeq
The least-favorable covariance matrix $\tilde{K}_{x_{t+1}}$ can be evaluated as
follows. Let
\bea
P_{t+1} &=& K_{x_{t+1}} - K_{x_{t+1} y_t} K_{y_t}^{-1} K_{y_t x_{t+1}} \nnl
    &=& (A_t -G_t C_t) V_t (A_t - G_t C_t)^T + (B_t -G_t D_t)(B_t - G_t D_t)^T 
\label{4.28}
\eea
and
\beq
V_{t+1} = \tilde{K}_{x_{t+1}} - K_{x_{t+1} y_t} K_{y_t}^{-1} K_{y_t x_{t+1}} 
\label{4.29}
\eeq
denote the nominal and least-favorable conditional covariance matrices of $x_{t+1}$
given ${\cal Y}_t$. Then 
\beq
V_{t+1}^{-1} = P_{t+1}^{-1} - \lambda_t^{-1} I_n \: , \label{4.30}
\eeq
where the Lagrange multiplier $\lambda_t > r(P_{t+1})$ is selected such that
\beq
\gamma_t (\lambda_t ) = \frac{1}{2} \Big[ \mbox{tr} \, ((I_n -\lambda_t^{-1}
P_{t+1})^{-1} - I_n ) + \ln \det (I_n - \lambda_t^{-1} P_{t+1} ) \Big]  = c_t \: .
\label{4.31}
\eeq
where as indicated in (\ref{2.17}), $\gamma_t (\lambda_t )$ is monotone decreasing
over $(r(P_{t+1}),\infty)$ and has for range ${\mathbb R}^+$. Thus for any
divergence tolerance $c_t >0$, there exists a matching Lagrange multiplier 
$\lambda_t > r(P_{t+1})$.     
\vskip 2ex

\noindent
{\bf Summary:} The least-favorable conditional distribution of $x_{t+1}$ given
$Y_t$ is given by
\beq
\tilde{f}_{t+1} (x_{t+1}|Y_t) \sim {\cal N} (\hat{x}_{t+1},V_{t+1})
\: , \label{4.32}
\eeq
where the estimate $\hat{x}_{t+1}$ is obtained by propagating the filter
(\ref{4.26})--(\ref{4.27}) and the conditional covariance matrix $V_{t+1}$
is obtained from (\ref{4.28}) and (\ref{4.30}), with the Lagrange multiplier
$\lambda_t$ specified by (\ref{4.31}). By writing $\theta_t = \lambda_t^{-1}$,
we recognize immediately that the robust filter is a form of risk-sensitive
filter of the type discussed in \cite{Whi} \cite[Chap. 10]{SC}. However, there 
is a new twist in the sense that, whereas standard risk-sensitive filtering
uses a fixed risk sensitivity parameter $\theta$, here $\theta$ is
time-dependent. Specifically, in classical risk-sensitive filters, $\theta$ 
is an exponentiation parameter appearing in the exponential of quadratic
cost to be minimized. Similarly, in earlier works \cite{BJP,LN,YUP,HS1,HS}     
relating risk sensitive filtering with minimax filtering with a relative
entropy constraint, a single global relative entropy constraint is
imposed, resulting in a single Lagrange multiplier/risk sensitivity
parameter. Here each component $\phi_t (z_t|x_t)$ of the model has an 
associated relative entropy constraint (\ref{3.10}), where the tolerance
$c_t$ varies in inverse proportion with the modeller's confidence in
the model component. In this respect, even if the state-space model
(\ref{3.1})--(\ref{3.2}) is time-invariant ($A$, $B$, $C$ and $D$ are
constant) and the tolerance $c_t = c$ is constant, the risk sensitivity
parameter $\lambda_t^{-1}$ will be generally time-varying. On the other hand,
if we insist on holding $\theta = \lambda^{-1}$ constant, it means that
the tolerance $c_t = \gamma_t (\lambda)$ is time varying since the
covariance matrix $P_{t+1}$ given by (\ref{4.28}) depends on time. 

\setcounter{equation}{0}
\section{Least-Favorable Model}
\label{sec:leastfav}

In last section, we derived the robust filter, which is of course
the most important component of the solution of the minimax filtering
problem. However for simulation and performance evaluation purposes, 
it is also useful to construct the least favorable model corresponding
to the optimum filter. Before proceeding, note that if $e_t = x_t -\hat{x}_t$
denotes the state estimation error, by subtracting (\ref{4.26})
from the state dynamics (\ref{3.1}) and taking into account
expression (\ref{3.2}) for the observations, the estimation 
error dynamics are given by
\beq
e_{t+1} = (A_t -G_t C_t) e_t + (B_t -G_t D_t) v_t \: , \label{5.1}
\eeq
where in the nominal model, the driving noise $v_t$ is independent
of error $e_t = x_t -\hat{x}_t$, since $\hat{x}_t$ depends exclusively
on observations $\{ y(s) ,~ 0 \leq s \leq t-1 \}$.

To find the least-favorable model, we use the characterization (\ref{4.4})
where $g_t = g_t^0$ is given by the robust filter (\ref{4.26}). This
gives
\beq
\tilde{\phi}_t^0 (z_t|x_t) = \frac{1}{M_t (\lambda_t)} \exp \Big( 
\frac{||e_{s+1}||^2}{2\lambda_t} \Big)
\phi_t (z_t|x_t) \: . \label{5.2}
\eeq
At this point, recall that $\tilde{\phi}_t^0 (z_t|x_t)$ is an unormalized
density. Specifically, integrating it over $z_t$ does not yield one,
but as we shall see below, a positive function of $e_t  =x_t -\hat{x}_t$. This
feature indicates that the maximizing player has the opportunity to change
retroactively the least-favorable density of $x_t$ (and therefore of earlier 
states) by selecting the model component $\tilde{\phi}_t^0(z_t|x_t)$. 
Properly accounting for this retroactive change forms an important aspect 
of the derivation of the least-favorable model. Instead of attempting 
to characterize directly the least-favorable density  of $z_t$ given
$x_t$, it is easier to find the least-favorable density of the
driving noise
\beq
v_t = \Gamma_t^{-1} \Big( z_t - \bmat {c} 
A_t \\ C_t 
\emat x_t \Big) \: . \label{5.3}
\eeq
Given $x_t$, the transformation (\ref{5.3}) establishes a one-to-one
correspondence between $z_t$ and $v_t$, so there is no loss of information
in characterizing the least-favorable model in terms of $v_t$. Let 
$\psi_t (v_t)$ and $\tilde{\psi}_t (v_t|e_t)$ denote respectively the 
nominal and least-favorable densities of noise $v_t$, where as will
be shown below, $\tilde{\psi}_t (v_t|e_t)$ actually depends on $e_t$. The 
nominal distribution is given by 
\beq
\psi_t (v_t) = \frac{1}{(2\pi)^{n+p}} \exp (-||v_t||^2 /2) \: .
\label{5.4}
\eeq

Assume that we seek to construct the least-favorable model of $z_t$
over a fixed interval $0 \leq t \leq T$, and that the least-favorable
noise distribution $\tilde{\psi}_s (v_s|e_s)$ has been identified
for $t+1 \leq s \leq T$. Accordingly, we have
\beq
\prod_{s=t+1}^T \exp \Big(\frac{||e_{s+1}||^2}{2\lambda_s} \Big) 
\psi_s (v_s) \sim \exp \big(||e_{t+1}||_{\Omega_{t+1}^{-1}}^2 /2 \big)
\prod_{s=t+1}^T \tilde{\psi}_s (v_s|e_s) \: , \label{5.5} 
\eeq
where $\sim$ indicates equality up to a multiplicative constant.
The term $\exp \big( ||e_{t+1}||_{\Omega_{t+1}^{-1}}^2 /2)$ appearing
in the above expression accounts for the cumulative effect of retroactive 
probability density changes performed by the maximizing player. Here $\Omega_t$ 
denotes a positive definite matrix of dimension $n$ which is evaluated 
recursively. Then the least-favorable model $\tilde{\psi}_t 
(v_t|e_t)$ is obtained by backward induction. Decrementing the index
$t$ by $1$ in (\ref{5.5}) gives the identity
\beq 
\exp \Big ( \big( ||e_{t+1}||_{\Omega_{t+1}^{-1}}^2 + \lambda_t^{-1} ||e_{t+1}||^2 
\big)/2\Big) \psi_t (v_t) \sim \tilde{\psi}_t (v_t|e_t) 
\exp \big( ||e_t||_{\Omega_t^{-1}}^2 /2 \big) \: . \label{5.6} 
\eeq
Let 
\beq
W_{t+1} \defn ( \Omega_{t+1}^{-1} + \lambda_t^{-1} I_n )^{-1} \: . \label{5.7}
\eeq 
Then by substituting the error dynamics (\ref{5.1}), the left hand side
of identity (\ref{5.6}) becomes
\beq
\exp \Big( \big( ||(A_t -G_tC_t)e_t + (B_t -G_t D_t)v_t)||_{W_{t+1}^{-1}}^2 - 
||v_t||^2 \big)/2 \Big) \: , \label{5.8}
\eeq
and the right-hand side of (\ref{5.6}) is obtained by decomposing
the quadratic exponent of (\ref{5.8}) as a sum of squares in $v_t$
and $e_t$. By doing so, we find that the least-favorable noise
density is given by 
\beq
\tilde{\psi}_t (v_t|e_t) \sim {\cal N} (H_t e_t, \tilde{K}_{v_t} )
\: , \label{5.9}
\eeq
where
\beq
\tilde{K}_{v_t} = \Big( I_{n+p} - (B_t -G_t D_t)^T W_{t+1}^{-1} (B_t -G_t D_t) \Big)^{-1}
\label{5.10}
\eeq
and
\beq
H_t = \tilde{K}_{v_t} (B_t-G_tD_t)^T W_{t+1}^{-1}(A_t-G_tC_t) \: .
\label{5.11}
\eeq 
Thus the least-favorable density of the noise $v_t$ involves a perturbation of
both the mean and the variance of the nominal noise distribution. The mean perturbation
is proportional to the filtering error $e_t$, which creates a coupling beetween
the robust filter and the least favorable model specified by dynamics and 
observations (\ref{3.1})--(\ref{3.2}) and least-favorable noise statistics
(\ref{5.9})--(\ref{5.11}).   

Finally, by matching quadratic components in $e_t$ on both sides of (\ref{5.6}),
we find
\beq
\Omega_t^{-1} = (A_t -G_t C_t)^T W_{t+1}^{-1} (A_t -G_t C_t) - H_t^T \tilde{K}_{v_t} H_t
\label{5.12}
\eeq
where $\tilde{K}_{v_t}$ and $H_t$ are given by (\ref{5.10}) and (\ref{5.11}).
By using the matrix inversion lemma \cite[p. 48]{Lau}
\[
(\alpha + \beta \gamma \delta)^{-1} = \alpha^{-1} -\alpha^{-1} \beta
(\gamma^{-1} + \delta \alpha^{-1}\beta )^{-1} \delta \alpha^{-1} 
\]
with $\alpha = W_t$, $\beta = (B_t -G_t D_t)$, $\gamma = -I_{n+p}$ and
$\delta = (B_t-G_t D_t)^T$ on the right-hand side of (\ref{5.12}), we obtain
\beq
\Omega_t^{-1} = (A_t -G_t C_t)^T [W_{t+1} -(B_t-G_t D_t)(B_t-G_t D_t)^T]^{-1} (A_t -G_t C_t)
\: . \label{5.13}
\eeq  
The recursion  (\ref{5.13}), together with (\ref{5.7}) specifies a backward
backward recursion which is used to account for retroactive changes of
previous least-favorable model densities performed by the maximizing player.  
The backards recursion is initialized with $\Omega_{T+1}^{-1} =0$, or 
equivalently,
\beq
W_{T+1} = \lambda_T I_{n+p} \: . \label{5.14}
\eeq
In this respect, it is interesting to note that recursion (\ref{5.13})
can be rewritten in the forward direction as 
\beq
W_{t+1} = (A_t -G_t C_t)\Omega_t(A_t-G_tC_t)^T + (B_t -G_t D_t)(B_t-G_tD_t)^T
\label{5.15}
\eeq
and (\ref{5.7}) is of course equivalent to
\beq
\Omega_{t+1} = (W_{t+1}^{-1} - \lambda_t^{-1} I_n )^{-1} \: .  
\label{5.16}
\eeq
Thus $\Omega_t$ and $W_t$ obey exactly the same forward recursions 
as as $V_t$ and $P_t$, but they are computed in the backward 
direction. Indeed, observe that the matrices $\Omega_t^{-1}$
and $W_t^{-1}$ are typically very small, so it is much easier to
maintain positive-definiteness by using (\ref{5.7}) which accumulates
small positive terms, instead of using (\ref{5.16}) which subtracts
a small positive-definite matrix from another one.
 
The least-favorable noise model (\ref{5.9})--(\ref{5.11}) can be viewed as a
time-varying version of the least-favorable model derived asymptotically
for the case of a constant model by Hansen and Sargent in \cite[Sec. 17.7]{HS}.
Specifically, the dynamics of the least-favorable model described in \cite{HS} 
are expressed in terms of the solution of a deterministic infinite-horizon
linear-quadratic regulator problem. The counterpart of this regulator is formed
here by backward recursion (\ref{5.13}), (\ref{5.7}). 

The model (\ref{5.9})--(\ref{5.11}) indicates that the driving noise $v_t$ 
admits the representation
\beq
v_t = H_t e_t +L_t \epsilon_t \label{5.17}
\eeq
where $L_t$ is an arbitrary matrix square root of $\tilde{K}_{v_t}$, i.e.,
\beq
L_t L_t^T = [ I_{n+p} - (B_t -G_t D_t)^T W_{t+1}^{-1} (B_t -G_t D_t)]^{-1}
\: , \label{5.18}
\eeq
and $\epsilon_t$ is a zero-mean WGN of variance $I_{n+p}$. Accordingly, as was
previously observed in \cite{HS}, if
\beq
\xi_t \defn \bmat {c}
x_t \\ e_t 
\emat \: , \label{5.19}
\eeq
the least-favorable model admits a state-space representation 
\bea
\xi_{t+1} &=& \tilde{A}_t \xi_t + \tilde{B}_t \epsilon_t \nnl
y_t &=& \tilde{C}_t \xi_t + \tilde{D}_t \epsilon_t \label{5.20}
\eea
with twice the dimension of the nominal state-space model, where
\bea
\tilde{A}_t = \bmat {cc} 
A_t & B_t H_t \\
0 & A_t -G_t C_t + (B_t -G_t D_t)H_t
\emat &,& 
\tilde{B}_t = \bmat  {c}
B_t  \\ 
(B_t - G_t D_t) 
\emat L_t \nnl
\tilde{C}_t = \bmat {cc}
C_t & D_tH_t 
\emat &,& 
\tilde{D}_t = D_t L_t \: .  \label{5.21}
\eea
Note that the model (\ref{5.20})--(\ref{5.21}) is constructed by 
performing first a forward sweep of the risk-sensitive filter
(\ref{4.27})--(\ref{4.28}) over interval $[0,T]$ to generate 
the gains $G_t$, followed by a backward sweep used to evaluate 
the matrix sequence $W_t$. Thus, the least-favorable model is 
constructed in a {\it nonsequential} manner, since increasing
the simulation interval beyond $[0,T]$ requires performing a new
backward sweep of recursion (\ref{5.13}), (\ref{5.7}). 

The model (\ref{5.20}) can be used to assess the performance of any
estimation filter designed under the assumption that the nominal model
(\ref{3.1})--(\ref{3.2}) is valid. Let $G_t^{\prime}$ be an arbitrary
time-dependent gain sequence, and let $\hat{x}_t^{\prime}$ be the state
estimate generated by the recursion
\beq
\hat{x}_{t+1}^{\prime} = A_t \hat{x}_t^{\prime} + G_t^{\prime} (y_t - C_t
\hat{x}_t^{\prime}) \: .
\label{5.22}
\eeq
Let $e_t^{\prime}  = x_t - \hat{x}_t^{\prime}$ denote the corresponding
filtering error. When the actual data is generated by the least-favorable
model (\ref{5.20})--(\ref{5.21}), by subtracting recursion (\ref{5.22})
from the first component of the state dynamics (\ref{5.20}), we obtain
\beq
\bmat {c}
e_{t+1}^{\prime} \\
e_{t+1}
\emat = \Big( \tilde{A}_t - \bmat {c}
G_t^{\prime} \\ 0
\emat \tilde{C}_t \Big) \bmat {c}
e_t^{\prime} \\ e_t
\emat + \Big( \tilde{B}_t - \bmat {c}
G_t^{\prime} \\ 0
\emat \tilde{D}_t \Big) \epsilon_t \: . \label{5.23}
\eeq
The recursion (\ref{5.23}) can be used to evaluate the performance
of filter (\ref{5.22}) when the data is generated by the least-favorable
model (\ref{5.20})--(\ref{5.21}). Specifically, consider the covariance
matrix
\[
\Pi_t = \tilde{E} \Big[ \bmat {c}
e_t^{\prime} \\ e_t \emat
\bmat {cc}
(e_t^{\prime})^T e_t^T \emat \Big] \: .
\]
By using the dynamics (\ref{5.23}) derived under the assumption that
the data is generated by the least-favorable model, we obtain the
Lyapunov equation
\bea
\Pi_{t+1} &=& \Big (\tilde{A}_t - \bmat {c}
G_t^{\prime} \\ 0
\emat \tilde{C}_t \Big) \Pi_t \Big ( \tilde{A}_t - \bmat {c}
G_t^{\prime} \\ 0
\emat \tilde{C}_t \Big)^T \nnl
&& \hspace*{0.2in} + \Big( \tilde{B}_t - \bmat {c}
G_t^{\prime} \\ 0
\emat \tilde{D}_t \Big) \Big( \tilde{B}_t - \bmat {c}
G_t^{\prime} \\ 0
\emat \tilde{D}_t \Big)^T \: , \label{5.24}
\eea
which can be used to evaluate the performance of filter (\ref{5.23})
when it is applied to the least-favorable model. For the special case
where $G_t^{\prime}$ is the Kalman gain sequence, this yields the
performance of the standard Kalman filter.

\setcounter{equation}{0}
\section{Simulations}
\label{sec:simu}

To illustrate the behavior of the robust filtering algorithm specified
by (\ref{4.26})--(\ref{4.31}), we consider a constant state-space model 
employed earlier in \cite{PS,Say}:
\bea
A = \bmat {cc} 
0.9802 & 0.0196 \\
0 & 0.9802
\emat &,& BB^T = Q = \bmat {cc}
1.9608 & 0.0195 \\
0.0195 & 1.9605
\emat \nnl
C = \bmat {cc} 1 & -1 \emat
&,& DD^T = 1 \: . \nn
\eea
The nominal process noise $Bv_t$ and measurement noise $Dv_t$ are assumed 
to be uncorrelated, so that $B D^T =0$, and the initial value of the
least-favorable error covariance matrix is selected as 
\[
V_0 = I_2 \: .
\]
We apply the robust filtering algorithm over an interval of length
$T=200$ for progressively tighter values $10^{-2}$, $10^{-3}$ and $10^{-4}$ 
of the relative entropy tolerance $c$. The corresponding time-varying risk-sensitivity
parameters $\theta_t = \lambda_t^{-1}$ obtained from (\ref{4.31}) are plotted
in \fig{thetol}. The least-favorable variances (the (1,1) and (2,2) entries of
$V_t$) of the two states are plotted as functions of time in \fig{x1tol}  
and \fig{x2tol}. As can be seen from the plots, although the relative 
entropy tolerance bounds that we consider are small, increasing the 
tolerance $c$ by a factor $10$ leads to an increase of about 7dB in
the state error variances.

\bef[htb!]
\centering
\includegraphics[width = 4in, height =4in]{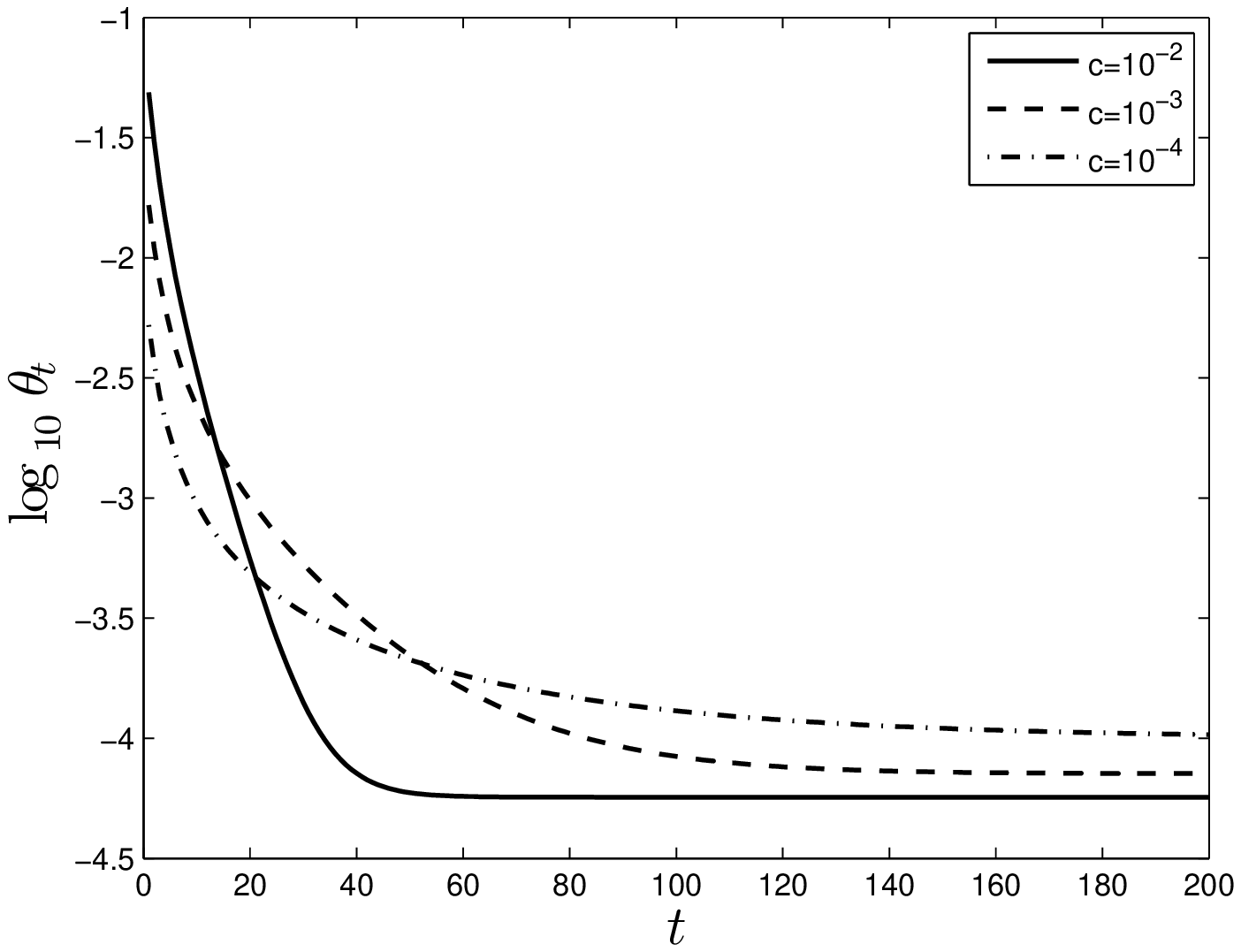}
\caption{Plot of time varying parameter $\theta_t = \lambda_t^{-1}$
(logarithmic scale) for $c=10^{-2}$, $10^{-3}$ and $10^{-4}$.}
\label{thetol}
\eef

\bef[htb!]
\centering
\includegraphics[width = 4in, height =4in]{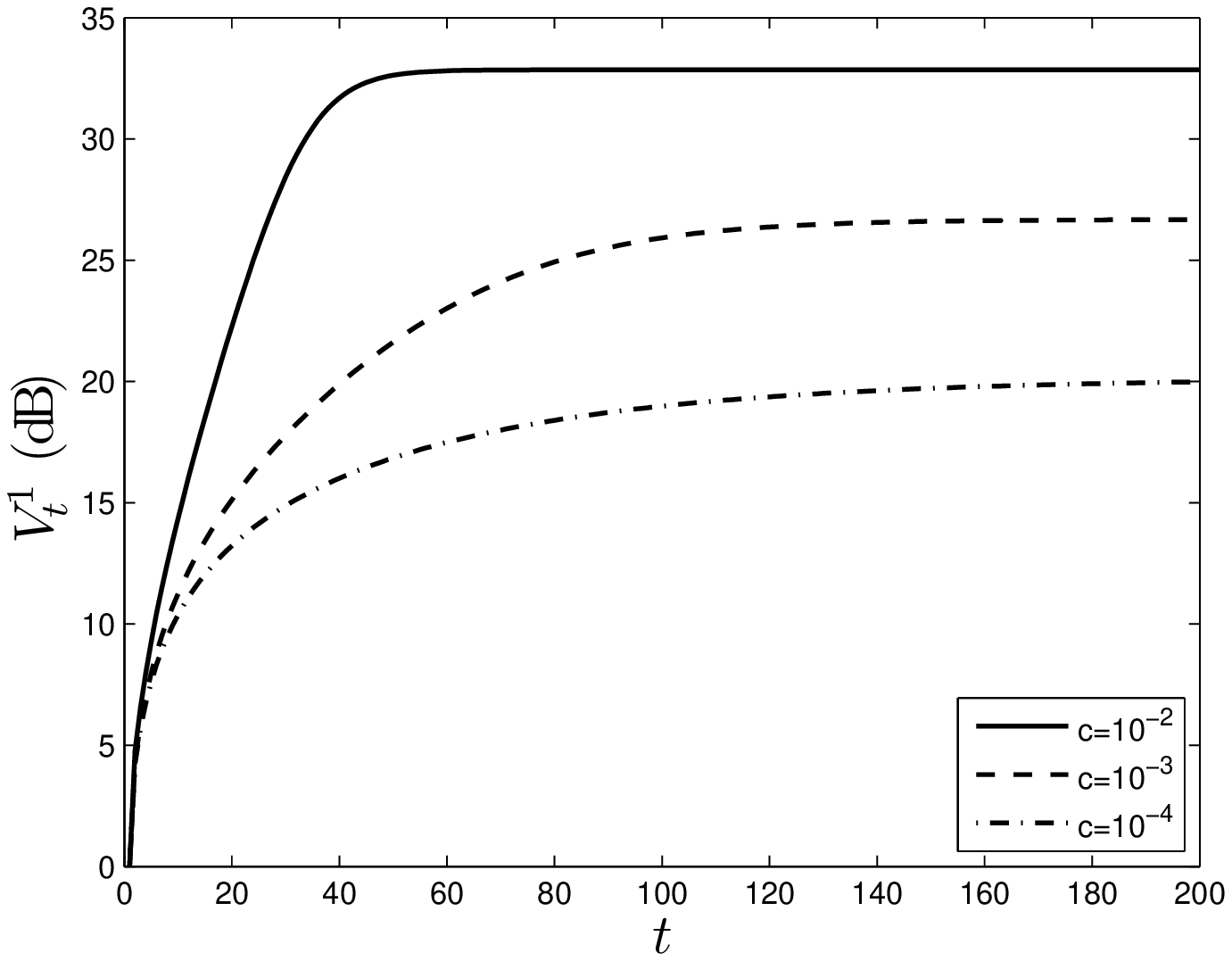}
\caption{Error variance of $x_{1t}$ (dB scale) for $c=10^{-2}$, 
$10^{-3}$, and $10^{-4}$.}
\label{x1tol}
\eef

\bef[htb!]
\centering
\includegraphics[width = 4in, height =4in]{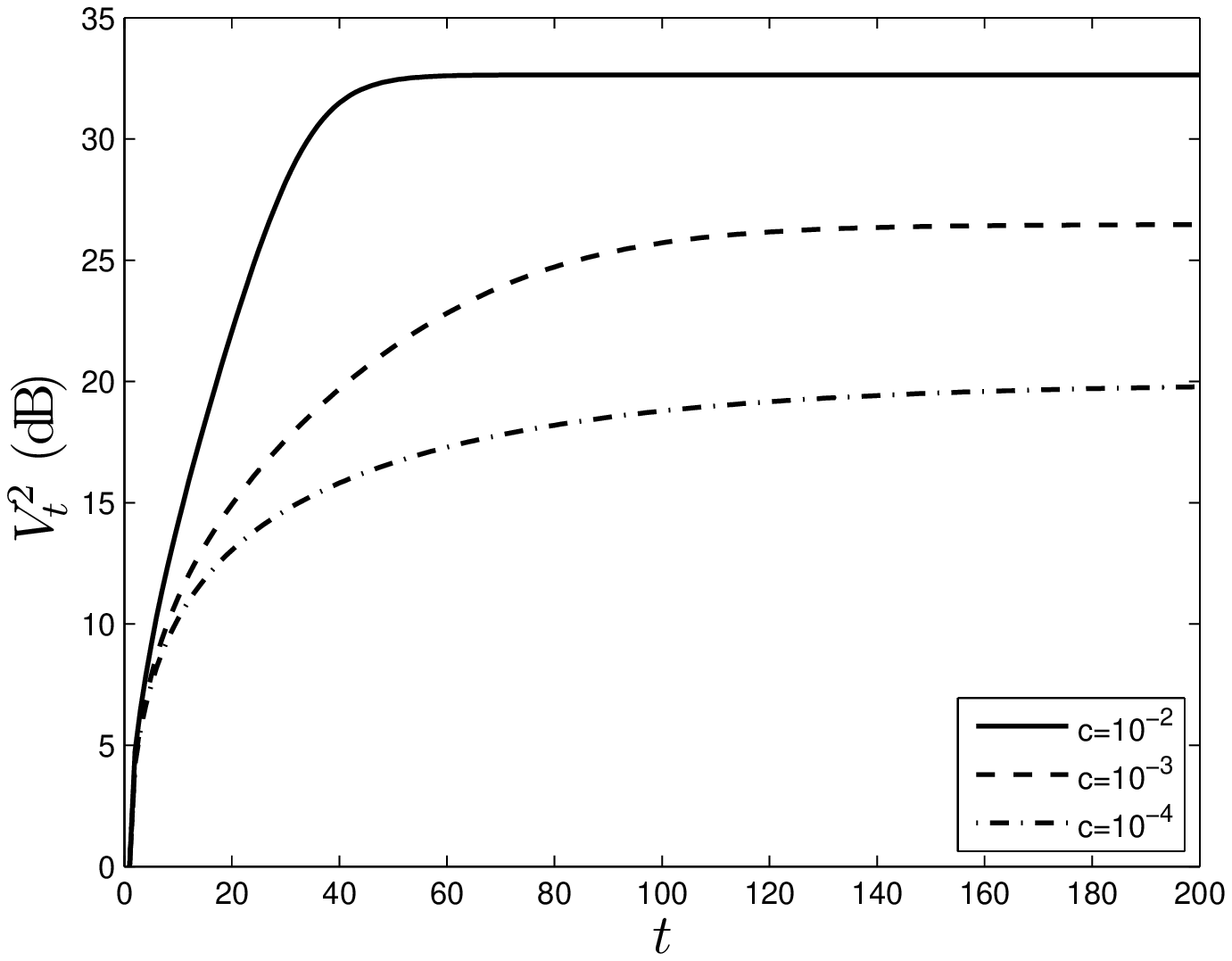}
\caption{Error variance of $x_{2t}$ (dB scale) for $c=10^{-2}$, 
$10^{-3}$, and $10^{-4}$.}
\label{x2tol}
\eef

Next, for a tolerance $c=10^{-4}$, we compare the performance of 
the risk-sensitive and Kalman filters for the nominal model,
and for the least-favorable model constructed as indicated in 
Section \ref{sec:leastfav}. The variances of the two-states
for the nominal model are shown in \fig{nomi1} and \fig{nomi2},
respectively. Clearly, the loss of performance of the 
risk-sensitive filter compared to the Kalman filter is less
than 1dB.

\bef[htb!]
\centering
\includegraphics[width = 4in, height =4in]{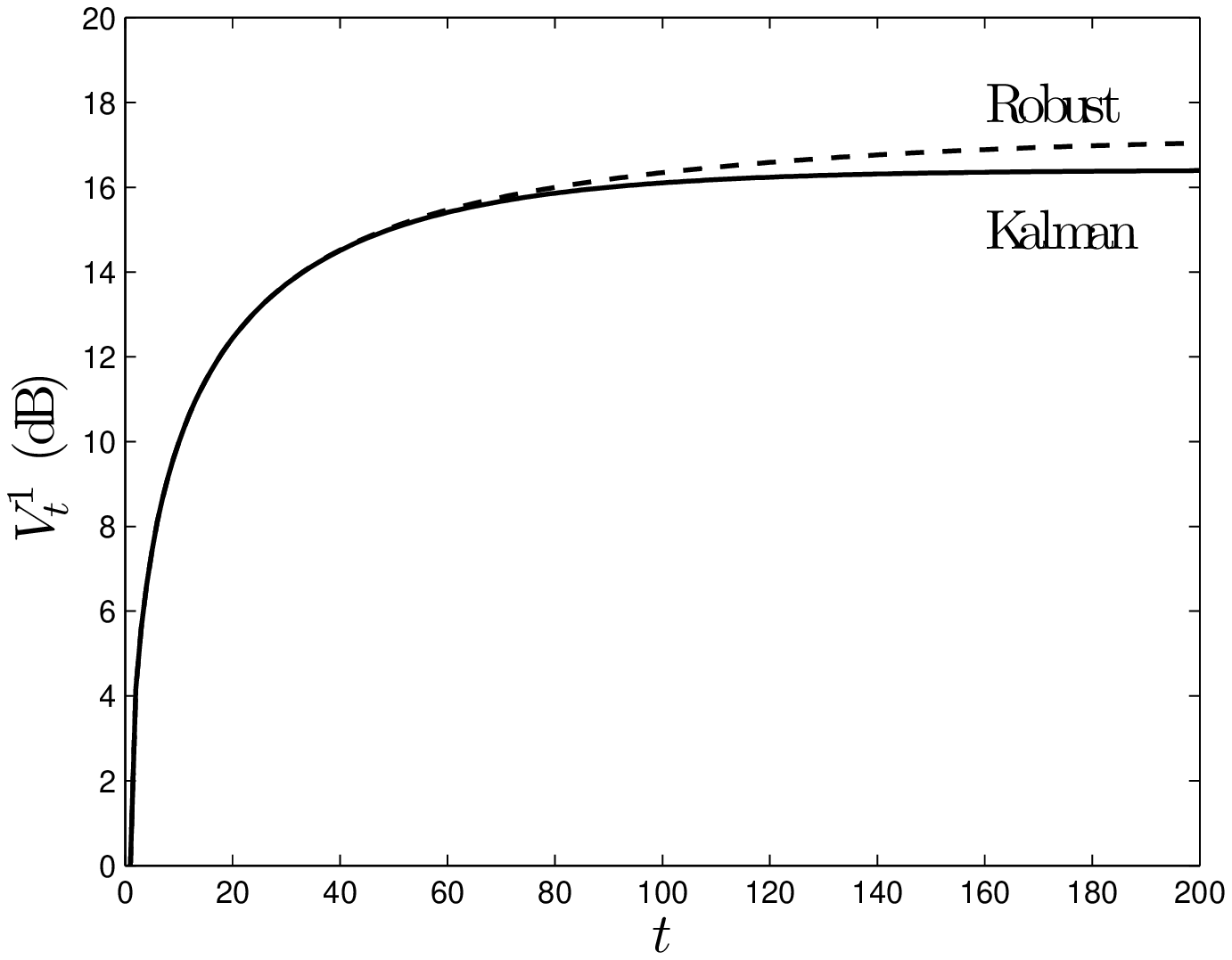}
\caption{Error variance of $x_{1t}$ (dB scale) when the 
risk-sensitive filter with $c=10^{-4}$ and the Kalman filter
are applied to the nominal model.} 
\label{nomi1}
\eef

\bef[htb!]
\centering
\includegraphics[width = 4in, height =4in]{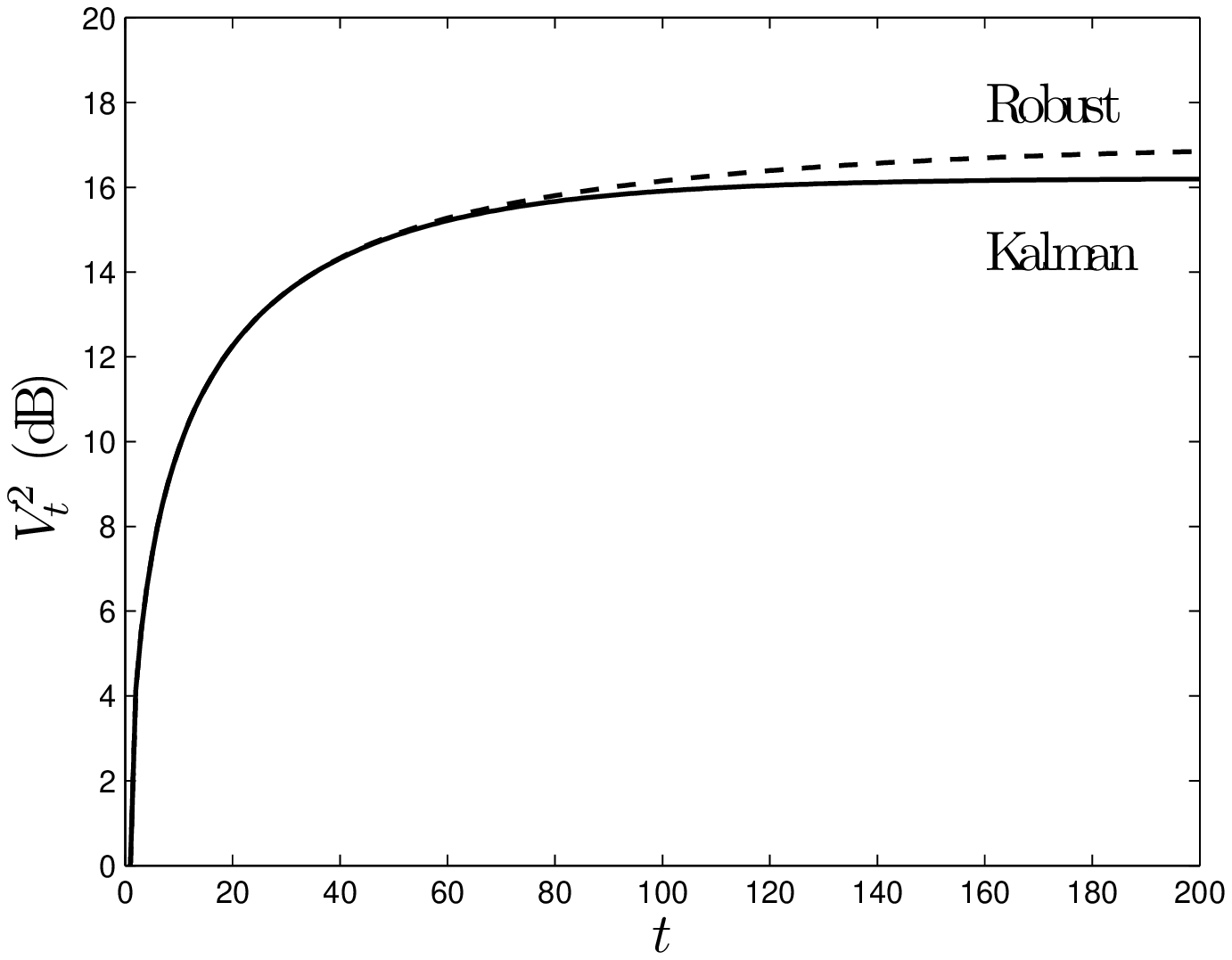}
\caption{Error variance of $x_{2t}$ (dB scale) when the 
risk-sensitive filter with $c=10^{-4}$ and the Kalman filter
are applied to the nominal model.} 
\label{nomi2}
\eef

On the other hand, as indicated in \fig{leastfav1} and \fig{leastfav2},
when the risk-sensitive and Kalman filters are applied to the
least-favorable model, the Kalman filter performance is about
8dB worse than the robust filter. Note that to allow the
backward recursion (\ref{5.13}), (\ref{5.7}) to reach steady
state, the backward model is computed for a larger interval, 
and only the first 200 samples of the simulation interval 
are retained, since later samples are affected by transients of 
the least-favorable model.

\bef[htb!]
\centering
\includegraphics[width = 4in, height =4in]{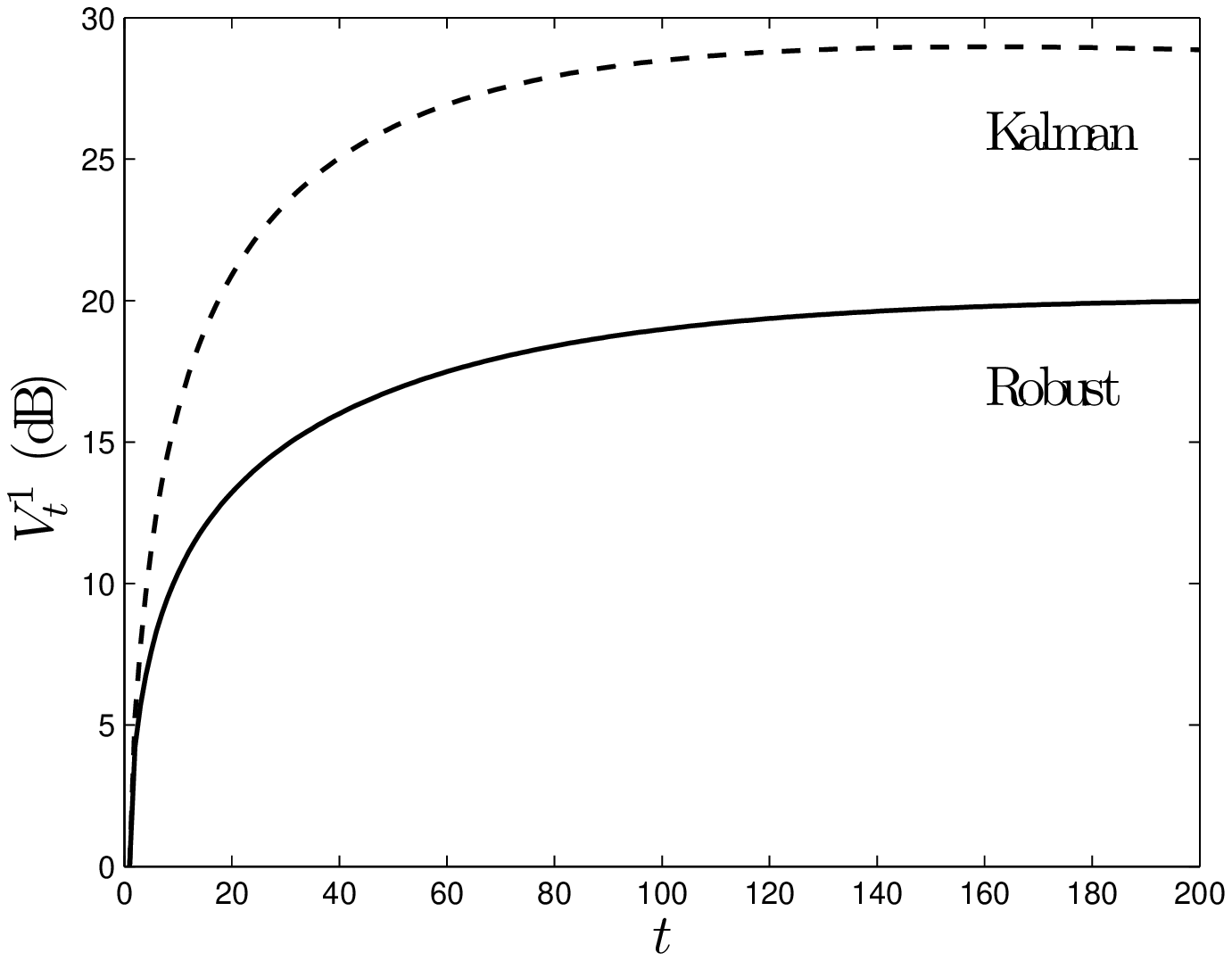}
\caption{Error variance of $x_{1t}$ (dB scale) when the 
risk-sensitive filter with $c=10^{-4}$ and the Kalman filter
are applied to the least-favorable model.} 
\label{leastfav1}
\eef

\bef[htb!]
\centering
\includegraphics[width = 4in, height =4in]{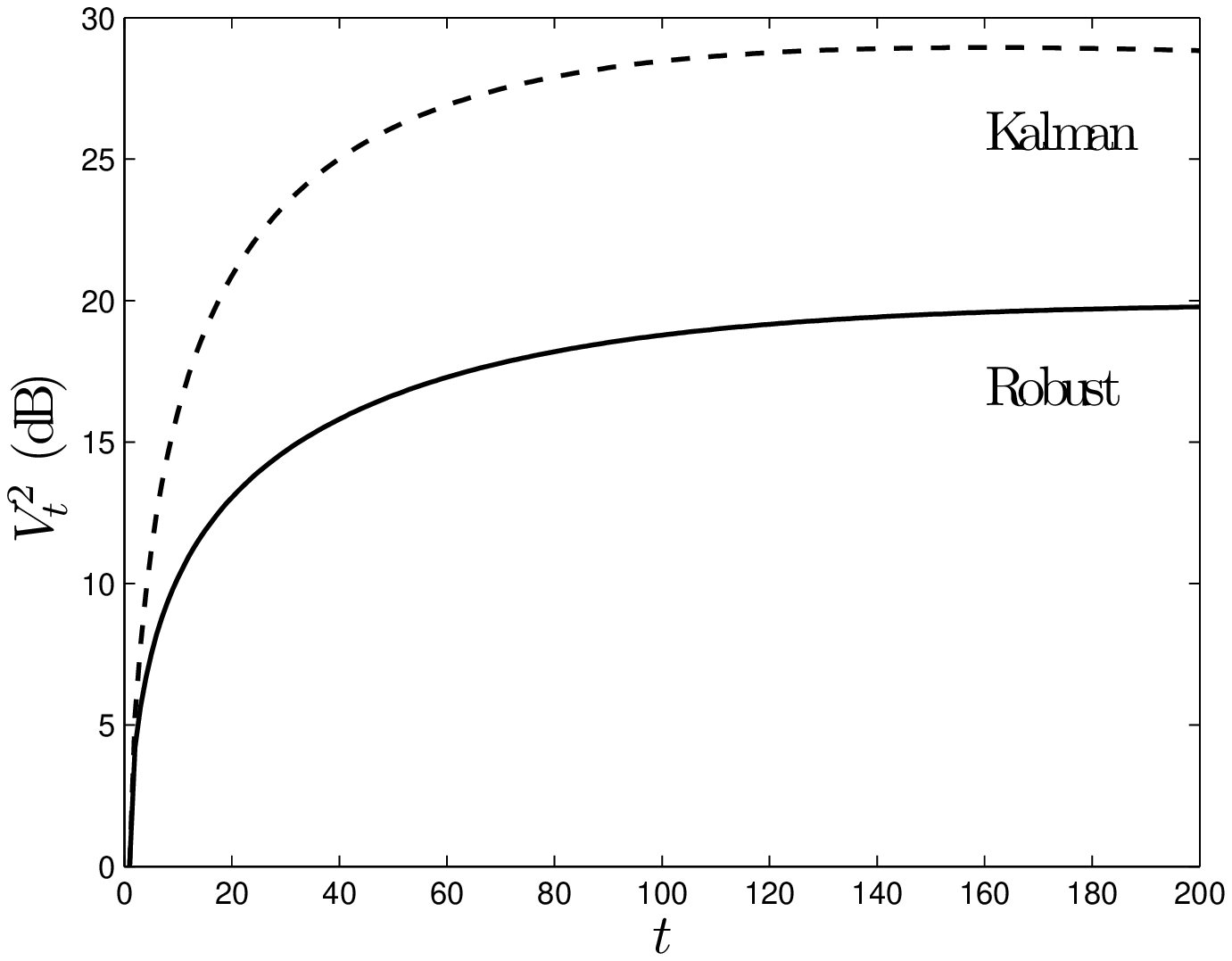}
\caption{Error variance of $x_{2t}$ (dB scale) when the 
risk-sensitive filter with $c=10^{-4}$ and the Kalman filter
are applied to the least-favorable model.} 
\label{leastfav2}
\eef

\section{Conclusion}
\label{sec:conc}

In this paper, we have considered a robust state-space filtering
problem with an incremental relative entropy constraint. The
problem was formulated as a dynamic minimax game, and by 
extending results presented in \cite{LN}, it was shown that
the minimax filter is a risk-sensitive filter with a time
varying risk-sensitive parameter. The associated least-favorable
model was constructed by performing a backward recursion which
keeps track of retroactive probability changes made by the maximizing 
player. The results obtained are similar in nature to those 
derived by Hansen and Sargent \cite{HS,HS1} for the minimax problem
(\ref{3.25}) when a single relative entropy constraint is applied 
to the overall state-space model and the maximizing agent is required
to operate under commitment..

A number of issues remain to be resolved. For the case of a constant
state-space model, it would be of interest to establish the 
convergence under appropriate conditions of the robust filtering 
recursions and of the backwards least-favorable model recursions.
One also has to wonder if the results derived here for
Gauss-Markov models could be extended to classes of systems, such 
as partially observed Markov chains, for which robust filtering
with an overall relative entropy constraint was considered 
previously in \cite{XUP}.   

\bibliography{robstate}

\end{document}